\documentclass[12pt,reqno]{amsart}
\usepackage{amsmath,amsfonts,amssymb,amscd,amsthm,amsbsy, color}
\usepackage[dvips]{graphicx}
\usepackage{comment}

\usepackage{wrapfig}
\usepackage{fancybox}

\usepackage{multimedia}
\usepackage{graphicx,epsfig}
\usepackage{ulem}

\numberwithin{equation}{section}

\newtheorem{theorem}{Theorem}
\newtheorem{meta-thm}[theorem]{Meta-Theorem}

\newtheorem{remark}[theorem]{Remark}
\newtheorem{definition}[theorem]{Definition}

\newcommand\beq[1]{ \begin{equation}\label{#1} }
\newcommand{\eeq}{ \end{equation} }

\newcommand\beqa[1]{ \begin{eqnarray} \label{#1}}
\newcommand{\eeqa}{ \end{eqnarray} }
\newcommand{\beqano}{ \begin{eqnarray*} }
\newcommand{\eeqano}{ \end{eqnarray*} }
\newcommand\equ[1]{{\rm (\ref{#1})}}

\def\integer{{\mathbb Z}}

\begin{document}

\title[On the dynamics of space debris: 1:1 and 2:1 resonances]
{On the dynamics of space debris: 1:1 and 2:1 resonances}

\author[A. Celletti]{Alessandra Celletti}

\address{
Department of Mathematics, University of Roma Tor Vergata, Via della Ricerca Scientifica 1,
00133 Roma (Italy)}

\email{celletti@mat.uniroma2.it}

\author[C. Gale\c s]{C\u at\u alin Gale\c s}

\address{
Department of Mathematics, Al. I. Cuza University, Bd. Carol I 11,
700506 Iasi (Romania)}
\email{cgales@uaic.ro}

\thanks{A.C. was partially supported by PRIN-MIUR 2010JJ4KPA$\_$009, GNFM-INdAM and by the European Grant MC-ITN Stardust.
}


\baselineskip=18pt              




\begin{abstract}
We study the dynamics of the space debris in the 1:1 and 2:1 resonances,
where geosynchronous and GPS satellites are located.
By using Hamiltonian formalism, we consider a model including the geopotential contribution for which we
compute the secular and resonant expansions of the Hamiltonian.

Within such model we are able to detect the equilibria and to study the main
features of the resonances in a very effective way.
In particular, we analyze the regular and chaotic behavior of the 1:1 and 2:1 resonant regions by analytical methods and
by computing the Fast Lyapunov Indicators, which provide a cartography of the resonances. This approach allows us to detect easily
the location of the equilibria, the amplitudes of the libration islands and the main dynamical stability features of the
resonances, thus providing an overview of the 1:1 and 2:1 resonant domains under the effect of Earth's oblateness.

The results are validated by a comparison with a model developed in Cartesian
coordinates, including the geopotential, the gravitational attraction of Sun and Moon and the solar radiation pressure.
\end{abstract}


\maketitle


\section{Introduction}\label{intro}
Since humankind started the conquest of space, a number of debris accumulated and now populate the
sky, from the near atmosphere to the geostationary region. Such debris, whose size runs from a few millimeters
to several centimeters, are remnants of spacecraft explosions or collisions, non-operative satellites, rocket upper stages, etc. (\cite{Klinkrad}).
Current estimates show that there exist about $3\cdot 10^8$ objects with size larger than 1 mm and
about $3.5\cdot 10^5$ objects larger than 1 cm. The impact of such debris with operative spacecraft or satellites
could result in a dangerous or even dramatic situation; the accumulation of debris in specific regions of the sky - where
most of operative satellites are positioned - cannot be neglected anymore. Understanding the dynamics and
evolution of space debris is essential for maintenance and control strategies, as well as to assess mitigation
procedures (\cite{RV2006}, \cite{RVF2000}).

The aim of this work is to provide a detailed study of a model describing the main resonant regions, precisely
where GEO and GPS satellites\footnote{GEO stands for \sl geostationary \rm orbit, located at about 42\,164 km from Earth's center,
while GPS is the acronym for \sl Global Positioning System, \rm a network of satellites at about 26\,560 km from Earth's center.} are positioned.
Our models are valid for spherical objects, typically space debris, while the study of satellites requires the analysis
of more complex effects, like the shape or flexibility of the spacecraft.
In particular, we introduce the Hamiltonian function describing
the effect of the Earth's gravity field and we expand it to get the secular and resonant terms.

Our main goal is to obtain an accurate, though computationally
fast, description of the resonant regions, using an approximation
of the 1:1 and 2:1 resonant Hamiltonians, so to be able to discern
between regular and chaotic behaviors (compare with
\cite{ElyHowell}, \cite{VL}, \cite{VLA}). We stress that the
mathematical tools adopted in the present work, although being
elementary, allow us to reach a twofold aim: to get a deeper
insight in the dynamics of resonant motions of space debris and,
at the same time, to obtain the main features of the dynamical
behavior in a reasonable computational time, especially when
compared to dynamical studies based on Cartesian equations. We
claim that the results provided by the Hamiltonian formalism
describing the geopotential effects can be used as starting point
to get a more detailed description of the dynamics. Furthermore,
we stress that the study presented in this work turns out to be
very effective also for the investigation of minor resonances
(e.g. 3:1, 3:2, 4:1, etc.), whose size is typically very small,
sometimes amounting to a few hundreds meters (compare with
\cite{CGminor}) and therefore extremely difficult to find with the
Cartesian formalism.

To describe the regions that will be the object of investigation of our study, let us introduce the following
(nowadays standard) classification of the sky, according to the altitude from the Earth,
in three main regions, referred to as LEO (acronym of Low--Earth--Orbit), MEO (Medium--Earth--Orbit) and GEO (Geostationary--Earth--Orbit); the altitude of such regions varies according to the following description.
LEO is the region between 0 and 2\,000 km, which is mainly affected by the Earth's attraction (including the effect of the oblateness)
and by the air drag. In this work we shall analyze objects in MEO (the region between 2\,000 and 30\,000 km) and GEO
(at altitudes larger than 30\,000 km, where the decay mechanisms are absent). Objects in MEO are mainly affected (in order of importance) by the $G M_E$ term, Earth's oblateness $J_2$ and $J_{22}$, the attraction of the Moon, the influence of
the Sun and the solar radiation pressure (hereafter SRP; see, e.g.,
\cite{chaogick}, \cite{deleflie2011}, \cite{rossi2008}, \cite{Sampaio}). In GEO the effects of Sun and Moon are
bigger than Earth's oblateness $J_{22}$ (see \cite{LDV}, \cite{VL}, \cite{VLA}, \cite{VLD}).
In the case of a large area--to--mass ratio, the solar radiation pressure is important and induces long periodic (yearly) variations in the eccentricity.

GPS orbits are located at about 26\,560 km from Earth's center; here satellites move with a
period of $12^h$ (sidereal time), namely half of the period of rotation of the Earth. We shall refer
to this situation as a \sl 2:1 gravitational resonance. \rm Much of the present work is devoted to the study of the 2:1 resonance,
whose typical effect is to provoke long--period changes of the eccentricity.

A geostationary orbit is situated on the equatorial plane at about 42\,164 km from Earth's center. Here, an object has an orbital period equal to the Earth's rotational period (one sidereal day), and thus it appears motionless, at a fixed position in the sky, to ground observers. This situation corresponds to a \sl 1:1 gravitational resonance. \rm

\vskip.1in

The 1:1 (in GEO) and 2:1 (in MEO) resonances will be the main object of this work, as, outside LEO, most of the satellites
(and therefore most of the debris) are currently found in these regions. Sometimes we will refer to these resonances
as GEO and MEO resonances, with clear reference to the satellites that populate the corresponding resonant regions.
The identification of the dynamical character around the 1:1 and 2:1
resonances is of seminal importance for the following reasons: stable regions allow us to minimize the eccentricity growth,
while chaotic behaviors can be skilfully exploited for mitigation strategies, either moving the debris in graveyard regions or
aiming at its re--enter and destruction into the atmosphere.

Most of the works available in the literature use Cartesian equations to describe the dynamics of space debris. The advantage of the Cartesian approach (see Section~\ref{sec:cart}) is that
one can easily include all effects (geopotential, lunisolar perturbations and solar radiation pressure; for other important
effects, like Earth's shadowing, see \cite{HL}). The
disadvantage is that it is difficult to catch the resonances, due to their limited size (often just of a few kilometers),
and to get a description of the dynamics inside the librational regions.
The present work aims to exploit the full power of the Hamiltonian formalism, which provides a careful and
detailed description of the resonant regions, once the main harmonics are taken into account. Since the Hamiltonian
approach contains some approximations, the results are validated by a comparison with the integration of the
complete Cartesian equations of motion.

Within the Hamiltonian formalism, we compute the secular and
resonant parts (see Section~\ref{sec:ham}), considering only the effect of Earth's oblateness
and disregarding all other effects. For these terms we provide explicit (sometimes long) expressions, which can be handled analytically or integrated numerically. While providing such expansions, we analyze the \sl dominant \rm terms
(namely those terms with magnitude bigger than other terms of the expansion) as a function of the eccentricity and the inclination. This procedure provides a very efficient
way to recognize which terms prevail in specific regions in the orbital parameters. The outcoming analysis shows a
marked difference between the 1:1 and 2:1 resonances, since for the 1:1 resonance a term of order $J_{22}$ dominates
almost everywhere, while for the 2:1 resonance different terms are dominant according to the value of the eccentricity
and the inclination. An expansion of the Hamiltonian around the resonant location and simple canonical
transformations allow us to reduce the Hamiltonian to a pendulum-like system with one single harmonic, which provides the
computation of the amplitude of the libration region around the resonance (see Section~\ref{sec:amplitude}). Though yielding an elementary estimate of the libration amplitudes, the results are in good agreement with the integrations performed using the complete Cartesian equations
(i.e. including the luni--solar attraction and the solar radiation pressure).

A cartography of the 1:1 and 2:1 resonances is provided in Section~\ref{sec:cartography} through the
computation of the Fast Lyapunov Indicators (hereafter FLIs). The FLIs provide an efficient tool to study the stable and chaotic behavior of a dynamical system by comparing the values of the FLIs as the initial conditions or parameters
are varied. Such analysis has been shown to be very effective in distinguishing between regular, resonant or chaotic motions
(\cite{froes}, \cite{froeslega}, see also  \cite{Guzzo1}, \cite{Guzzo2}, \cite{Gales}). Using the FLIs, different cartographic studies are performed in Section~\ref{sec:cartography} to describe
the 1:1 and 2:1 resonances (see also \cite{VDLC} where a different chaos indicator has been used).

For the 1:1 resonance we investigate the effects of Earth's gravitational perturbations up to degree and order four. In particular, we see that the location of the equilibrium points is not influenced by the longitude of the ascending node, but it depends on the argument of perigee. However, since for the 1:1 resonance a term of order $J_{22}$ dominates almost everywhere, this dependence on the argument of perigee is negligible, except for some specific regions in eccentricity and inclination, where more dominant terms coexist.

As far as the 2:1 resonance is concerned, by using a toy model including the secular term and a reduced number
of dominant terms, we show that there is a \sl superposition of harmonics. \rm   For inclinations different than the critical
inclination $i=63.4^o$ (at which so--called \sl frozen orbits \rm are found), we observe that the resonance splits
into a triplet of resonances with a complex interplay of regular and chaotic motions due to the interaction between
the different harmonics. It is also remarkable that, for
$i=70.53^o$, a \sl transcritical  bifurcation \rm phenomenon takes place: the stability of the equilibria reverses for small changes in the inclination, thus allowing to
move a debris from a stable to an unstable zone (or viceversa).
An exhaustive analysis of transcritical bifurcations for different resonances is performed in \cite{CGminor}.

This work is organized as follows. In Section~\ref{sec:cart} we introduce the Cartesian equations of
motion including the effects of the Earth's oblateness, the gravitational attraction of Sun and Moon,
and the solar radiation pressure (see also Appendix~\ref{app:geo}). The Hamiltonian approach is presented in Section~\ref{sec:ham}, which provides also the analysis of the resonant and dominant terms. The reduction to a pendulum--like structure is given in Section~\ref{sec:amplitude}, where an estimate of the resonant amplitudes is provided.
The cartographic study based on the computation of the FLIs is presented in Section~\ref{sec:cartography}.

\vskip.1in

\bf Acknowledgements. \rm We are grateful to Alessandro Rossi for very useful discussions and his constant encouragement;
we also thank Christoph Lhotka for helpful suggestions. We thank the reviewers for several comments which helped to improve
this work.

\section{Geopotential, lunisolar and solar radiation pressure effects}\label{sec:cart}

In this section we introduce the Cartesian equations of motion of
a space debris $S$ under the influence of the Earth, taking into
account the oblateness and  the rotation of our planet, the effects of
the Moon and of the Sun as well as the solar radiation pressure.
We refer to \cite{Beutler}, \cite{chao} for further details.

The equations of motion are provided by the sum of the
contributions of the Earth's gravitational influence, computed over the whole volume $V_E$ of the Earth and including
the oblateness effect, the solar and lunar attractions, proportional to the masses $m_S$ and $m_M$ of Sun and Moon,
and the solar radiation pressure:
\beqa{eq1}
\ddot{\mathbf{r}}&=&-G\int_{V_E} \rho(\mathbf{r}_p)\ {{\mathbf{r}-\mathbf{r}_p}\over {|\mathbf{r}-\mathbf{r}_p|^3}}\ dV_E - Gm_S \Bigl({{\mathbf{r}-\mathbf{r}_S}\over {|\mathbf{r}-\mathbf{r}_S|^3}}+{\mathbf{r}_S\over {|\mathbf{r}_S|^3}}\Bigr)\nonumber\\
&-& Gm_M\Bigl({{\mathbf{r}-\mathbf{r}_M}\over
{|\mathbf{r}-\mathbf{r}_M|^3}}+{\mathbf{r}_M\over
{|\mathbf{r}_M|^3}}\Bigr)+C_r P_r a_S^2\ ({A\over m})\ {{\mathbf{r}-\mathbf{r}_S}\over {|\mathbf{r}-\mathbf{r}_S|^3}}\ , \eeqa where
$G$ is the gravitational constant,
$\rho(\mathbf{r}_p)$ is the density at some point $\mathbf{r}_p$
inside the Earth, $\mathbf{r}$, $\mathbf{r}_S$, $\mathbf{r}_M$ are the position vectors of the debris, the Sun and the Moon
with respect to the Earth's center, $C_r$ is the reflectivity coefficient, depending on the optical properties
of the space debris surface, $P_r$ is the radiation pressure for an object located at $a_S=1$ AU, $A\over m$ is
the area--to--mass ratio with $A$ being the cross--section of the space debris and $m$ its mass. The vectors
$\mathbf{r}_S$, $\mathbf{r}_M$  can be found in \cite{Beutler}, \cite{MG} as a function of time.

We introduce two reference frames centered in the Earth, the quasi--inertial frame with unit vectors $\{\mathbf{e}_1,\mathbf{e}_2,\mathbf{e}_3\}$ fixed
and the synodic frame with unit vectors $\{\mathbf{f}_1,\mathbf{f}_2,\mathbf{f}_3\}$ rotating with the same angular velocity of the Earth.
We choose the orientation of these vectors such that $\mathbf{e}_3=\mathbf{f}_3$ is perpendicular to the equatorial plane and $\mathbf{f}_1$, $\mathbf{f}_3$ are in the Greenwich meridian plane.

Denoting by $(r,\lambda,\phi)$ the spherical coordinates in the synodic frame, the
geopotential associated to the first term at the right hand side of \equ{eq1},
expanded in spherical harmonics, can be written as (\cite{Beutler}, \cite{Kaula}, \cite{MG})
$$
V(r,\phi,\lambda)={{GM_E}\over r}\
\sum_{n=0}^\infty \Big({R_E\over r}\Big)^n\ \sum_{m=0}^n P_n^m(\sin\phi)\
(C_{nm}\cos m\lambda+ S_{nm}\sin m\lambda)\ ,
$$
where $M_E$, $R_E$ are the mass and equatorial radius of the Earth,
the quantities $P_n^m$ are defined in terms of the Legendre polynomials:
$$
P_n(x)\equiv {1\over {2^n n!}}\ {{d^n}\over {dx^n}}\{(x^2-1)^n\}\ ,\qquad
P_n^m(x)\equiv (1-x^2)^{m\over 2}\ {{d^m}\over {dx^m}}\{P_n(x)\}\ ,
$$
while $C_{nm}$, $S_{nm}$ are the spherical harmonic
coefficients of the geopotential (\cite{Kaula}).

For reasons which will be clear in the following sections, in our computations we consider a model with Earth's gravity harmonics up to degree and order $n=m=3$. Since the corresponding system of equations involves long expressions, in the following we present a simplified system obtained by taking into account just the Earth's gravity harmonics up to degree and order two; for self--consistency we provide in the Appendix some details on the derivation of the equations of motion. The reader can easily use the results presented there in order to get the full system of the equations that we use in this paper.
Denoting by $(x,y,z)$ the coordinates in the quasi--inertial frame, the equations of motion, in which the Earth's  gravity
harmonics are considered up to $n=m=2$, are:
\beqano
\ddot x&=&-{{GM_Ex}\over r^3}+{{GM_ER_E^2}\over r^5}\{C_{20}({3\over
2}x-{{15}\over 2}{{xz^2}\over r^2})+6 C_S^- x+6 C_S^+ y\nonumber\\
&+&{{15 x}\over r^2}[C_S^-(y^2-x^2)-2xy C_S^+]\}-Gm_S\Bigl({{x-x_S}\over {|\mathbf{r}-\mathbf{r}_S|^3}}+{x_S\over r_S^3}\Bigr)\nonumber\\
&-&Gm_M \Bigl({{x-x_M}\over {|\mathbf{r}-\mathbf{r}_M|^3}}+{x_M\over r_M^3}\Bigr)+C_r P_r a_S^2\ ({A\over m})\ {{x-x_S}\over {|\mathbf{r}-\mathbf{r}_S|^3}}
\eeqano
\beqano
\ddot y&=&-{{GM_Ey}\over r^3}+{{GM_ER_E^2}\over
r^5}\{C_{20}({3\over 2}y-{{15}\over 2}{{yz^2}\over r^2})+6 C_S^+
x-6 C_S^- y\nonumber\\
&+&{{15 y}\over r^2}[C_S^-(y^2-x^2)-2xy C_S^+]\}-Gm_S\Bigl({{y-y_S}\over {|\mathbf{r}-\mathbf{r}_S|^3}}+{y_S\over r_S^3}\Bigr)\nonumber\\
&-&Gm_M \Bigl({{y-y_M}\over {|\mathbf{r}-\mathbf{r}_M|^3}}+{y_M\over r_M^3}\Bigr)+C_r P_r a_S^2\ ({A\over m})\ {{y-y_S}\over {|\mathbf{r}-\mathbf{r}_S|^3}}
\eeqano
\beqa{eq3}
\ddot z&=&-{{GM_Ez}\over r^3}+{{GM_ER_E^2}\over
r^5}\{C_{20}({9\over 2}z-{{15}\over 2}{{z^3}\over r^2})
+{{15 z}\over r^2}[C_S^-(y^2-x^2)-2xy C_S^+]\}\nonumber\\
&-&Gm_S\Bigl({{z-z_S}\over
{|\mathbf{r}-\mathbf{r}_S|^3}}+{z_S\over r_S^3}\Bigr)-Gm_M
\Bigl({{z-z_M}\over {|\mathbf{r}-\mathbf{r}_M|^3}}+{z_M\over
r_M^3}\Bigr)\nonumber\\
&+&C_r P_r a_S^2\ ({A\over m})\ {{z-z_S}\over {|\mathbf{r}-\mathbf{r}_S|^3}}\ ,
\eeqa
where $C_S^-\equiv C_{22}\cos2\theta-S_{22}\sin2\theta$, $C_S^+\equiv C_{22}\sin 2\theta+S_{22}\cos 2\theta$,
and $\theta $ is the sidereal time.

\section{Hamiltonian formulation of the resonances under the geopotential}\label{sec:ham}
In this section we provide the Hamiltonian formulation of the equations of motion
by considering just the effect of the geopotential. We give explicit expressions for the
secular part as well as the 1:1 and 2:1 resonant Hamiltonians. We use action--angle Delaunay variables
$(L,G,H,M,\omega,\Omega)$, which are related to the orbital elements $(a,e,i,M,\omega,\Omega)$ by the expressions
\begin{equation}\label{Delaunay_var}
L=\sqrt{\mu_E a}\ , \qquad G=L \sqrt{1-e^2}\ , \qquad H=G \cos i\ ,
\end{equation}
where $\mu_E=GM_E$, $a$ is the semimajor axis, $e$ the eccentricity, $i$ the inclination,
$M$ the mean anomaly, $\omega$ the argument of perigee, $\Omega$ the longitude of the ascending node
(see, e.g., \cite{Alebook}).
Recalling that $\theta$ is the sidereal time, the Hamiltonian can be written as
\beq{H}
\mathcal{H}(L,G,H,M,\omega,\Omega,\theta)=-{\mu^2_E\over {2L^2}}+R_{earth}(a,e,i,M,\omega,\Omega,\theta)\ ,
\eeq
where $R_{earth}$ denotes the disturbing function, whose explicit expression will be given in Section~\ref{sec:geop}.

\subsection{The disturbing function $R_{earth}$}\label{sec:geop}
In the geocentric quasi--inertial frame, the disturbing function $R_{earth}$  is given by (see \cite{Kaula})
\beq{Rearth}
R_{earth}=- {{\mu_E}\over a}\ \sum_{n=2}^\infty \sum_{m=0}^n \Bigl({R_E\over a}\Bigr)^n\ \sum_{p=0}^n F_{nmp}(i)\
\sum_{q=-\infty}^\infty G_{npq}(e)\ S_{nmpq}(M,\omega,\Omega,\theta)\ .
\eeq
The functions $F_{nmp}$, $G_{npq}$ are given by the following relations (see \cite{Kaula}):
\beqa{Ffun}
F_{nmp}(i)&=&\sum_w {{(2n-2w)!}\over {w!(n-w)!(n-m-2w)!2^{2n-2w}}} \sin^{n-m-2w}i\ \sum_{s=0}^m\left(\begin{array}{c}
  m \\
  s \\
 \end{array}\right)
 \cos^si\nonumber\\
 &&\times \sum_c \left(\begin{array}{c}
  n-m-2w+s \\
  c \\
 \end{array}\right)
\left(\begin{array}{c}
  m-s \\
  p-w-c \\
 \end{array}\right)
 (-1)^{c-k}\ ,
\eeqa
where $k=[{{n-m}\over 2}]$, $w$ is summed from zero to the lesser of $p$ and $k$, $c$ is summed over all values for which the binomial coefficients are not zero; the functions $G_{npq}$ are defined as
\beq{Gfun}
G_{npq}(e)=(-1)^{|q|}(1+\beta^2)^n\beta^{|q|}\ \sum_{k=0}^\infty P_{npqk}Q_{npqk}\beta^{2k}\ ,
\eeq
where
\beqano
\beta&=&{e\over {1+\sqrt{1-e^2}}}\nonumber\\
P_{npqk}&=&\sum_{r=0}^h \left(\begin{array}{c}
  2p'-2n \\
  h-r \\
 \end{array}\right)
{{(-1)^r}\over r!}
({{(n-2p'+q')e}\over {2\beta}})^r\ ,
\eeqano
with $h=k+q'$ when $q'>0$ and $h=k$ when $q'<0$;
$$
Q_{npqk}=\sum_{r=0}^h \left(\begin{array}{c}
  -2p' \\
  h-r \\
 \end{array}\right)
{1\over r!}
({{(n-2p'+q')e}\over {2\beta}})^r\ ,
$$
where $h=k$ when $q'>0$ and $h=k-q'$ when $q'<0$,
$p'=p$ and $q'=q$ when $p\leq n/2$, $p'=n-p$ and $q'=-q$ when $p> n/2$.
It is worth mentioning that $G_{npq}(e)=\mathcal{O}(e^{|q|})$.

To complete the description of \equ{Rearth}, we provide the expression of $S_{nmpq}$.
If we introduce the quantities $J_{nm}$ and $\lambda_{nm}$ defined by
$$
J_{nm} = \sqrt{C_{nm}^2+S_{nm}^2}   \quad \textrm{if} \ m\neq 0\ , \qquad    J_{n0} \equiv J_n= -C_{n0}    \ ,$$
$$
 C_{nm}=-J_{nm} \cos(m \lambda_{nm}) \ , \qquad S_{nm}=-J_{nm} \sin(m \lambda_{nm}) \ ,
$$
then we can write $S_{nmpq}$ in the form
\beq{S}
S_{nmpq}=\left\{%
\begin{array}{cc}
  -J_{nm}  \cos \widetilde{\Psi}_{nmpq}  & \textrm{if} \ n-m\  \textrm{is even} \\
  -J_{nm}  \sin \widetilde{\Psi}_{nmpq}  & \ \textrm{if} \ n-m\ \textrm{is odd}\ , \\
 \end{array}%
\right.
\eeq
where
\beq{psitilde}
\widetilde{\Psi}_{nmpq}=(n-2p) \omega+(n-2p+q)M+m(\Omega-\theta)-m \lambda_{nm}\ .
\eeq

\subsection{Expansion of the Hamiltonian}\label{sec:exp}
With reference to \equ{H}, the long term variation of the orbital elements is governed by the secular and resonant terms
associated to the perturbing function $R_{earth}$, once we average over the non--resonant terms.

Therefore, in the rest of this section we focus our attention on the secular and resonant parts up to terms of degree and order $n = m = 4$. We will see in Section~\ref{sec:cartography} that the fourth degree provides a sufficiently accurate description
of the dynamics. Therefore, we approximate $R_{earth}$ by
$$R_{earth}=R^{sec}_{earth}+R_{earth}^{res}+R_{earth}^{nonres}\cong -\sum_{n=2}^4 \sum_{m=0}^4 V_{nm}\ ,$$
where $R^{sec}_{earth}$, $R_{earth}^{res}$, $R_{earth}^{nonres}$ are, respectively, the secular, resonant and non-resonant parts of the Earth's potential and
$$
V_{nm}=\frac{\mu_E R_E^n}{a^{n+1}} \sum_{p=0}^n F_{nmp}(i) \sum_{q=-\infty}^{\infty} G_{npq}(e) S_{nmpq}(M,\omega, \Omega, \theta)\ .
$$
In order to define the resonant contributions, we need the following definition.

\begin{definition}
A $p:q$ gravitational resonance for $p$, $q\in\integer\backslash\{0\}$ occurs when  the orbital period of the object and the period of Earth's rotation are commensurable in the ratio ${p\over q}$, namely the following relation holds:
\beq{res}
q\ \dot{M}-p\ \dot{\theta} = 0\ ,
\eeq
where $\dot M$ provides the mean motion of the object and $\dot{\theta}$ is the angular speed of the Earth's rotation.
\end{definition}

We stress that \equ{res} can be satisfied only approximately, namely within a specific accuracy, as it
happens, e.g., in spin--orbit or mean--motion resonances in Celestial Mechanics.

Since the frequencies $\dot{\omega}$, $\dot{\Omega}$ are small but not zero, then for a specific resonance the resonant angles $\widetilde{\Psi}_{nmpq}$ in \equ{psitilde} for different $n$, $m$, $p$, $q$ have zero derivative at different locations, thus
providing that each resonance splits into a multiplet of resonances. The exact location of the resonance for each component of the multiplet is obtained by using the exact relation $\dot{\widetilde{\Psi}}_{nmpq}=0$.

In the following sub-sections we provide the explicit expansions of the secular and 1:1 resonant parts of $R_{earth}$
up to the second order in the eccentricity, while the 2:1 resonant contribution is expanded up to the fourth order
in the eccentricity, due to the fact that many highly eccentric satellites are located in the 2:1 resonant region.
Though it would suffice to provide the formulae given in Sections~\ref{sec:geop}, \ref{sec:exp} to compute the secular and resonant
parts of $R_{earth}$, we believe worthwhile to give the explicit expressions, since their forms are seminal for the discussion of
the dominant terms as well as for the introduction of toy models, which will help to explain some features of the resonant dynamics.

\subsubsection{The secular part of the disturbing function $R_{earth}$}
Recalling the expression for $S_{nmpq}$ given in \equ{S}, the secular terms correspond to $m=0$ and $n-2p+q=0$.
Using \equ{Rearth}, the values for the functions $F$ and $G$ obtained from \equ{Ffun}, \equ{Gfun},
and the formula \equ{S} for $S_{nmpq}$, we get the following expression for the
secular part of the geopotential up to second order in the eccentricity:
\beqa{Rsec}
R_{earth}^{sec}&\cong&\frac{\mu_E R^2_E J_{2}}{a^3} \Bigl(\frac{3}{4} \sin^2 i -\frac{1}{2}\Bigr) (1-e^2)^{-3/2} \nonumber\\
&+&\frac{2\mu_E R^3_E J_{3}}{a^4} \Bigl(\frac{15}{16} \sin^3 i -\frac{3}{4} \sin i\Bigr) e (1-e^2)^{-5/2} \sin \omega \nonumber \\
&+&\frac{\mu_E R^4_E J_{4}}{a^5} \Bigl[\Bigl(-\frac{35}{32} \sin^4 i +\frac{15}{16} \sin^2 i\Bigr) \frac{3e^2}{2}(1-e^2)^{-7/2} \cos(2\omega) \nonumber \\
&+&
\Bigl(\frac{105}{64} \sin^4 i -\frac{15}{8} \sin^2 i+\frac{3}{8}\Bigr) (1+\frac{3e^2}{2})(1-e^2)^{-7/2} \Bigr]\ .
\eeqa
We remark that since $J_{2} \gg J_{3}$ and $J_{2} \gg J_{4}$, the $J_2$-term is dominant.

We report in Table~\ref{table:J} the values obtained according to the EGM2008 model of $C_{nm}$,
$S_{nm}$ and $J_{nm}$ in units of $10^{-6}$, as well as the values
of $\lambda_{nm}$ (\cite{EGM2008}, see also \cite{chao},
\cite{MG}).

\subsubsection{Resonance 1:1.} The 1:1 resonant terms correspond to $n-2p+q=m$, $m>0$. Retaining just the terms up to second order in eccentricity, we obtain:
\begin{equation}
\begin{split}
& R_{earth}^{res1:1}\cong \frac{\mu_E R_E^2 J_{22}}{a^3} \Bigl\{\frac{3}{4} (1+\cos i)^2 (1-\frac{5}{2}e^2) \cos [2(M-\theta +\omega +\Omega -\lambda_{22})]\nonumber\\
& \qquad + \frac{27}{8} e^2 \sin^2 i  \cos [2(M-\theta +\Omega -\lambda_{22})]\Bigr\}\nonumber\\
& \ + \frac{\mu_E R_E^2 J_{21}}{a^3} \Bigl\{ \frac{3}{4} \sin i (1+\cos i) (-\frac{e}{2}) \cos (M-\theta +2 \omega +\Omega -\lambda_{21}) \nonumber\\
& \qquad + \frac{3}{2} e (-{3\over 2}\sin i\cos i)  \cos [M-\theta +\Omega -\lambda_{21}]\Bigr\}\nonumber\\
& \ + \frac{\mu_E R_E^3 J_{31}}{a^4} \Bigl\{ -\frac{15}{16} \sin^2 i (1+\cos i) \frac{e^2}{8} \cos (M-\theta +3 \omega +\Omega -\lambda_{31}) \\
& \qquad + \Bigl(\frac{15}{16} \sin^2 i (1+3\cos i)-\frac{3}{4} (1+\cos i) \Bigr) (1+ 2 e^2) \cos (M-\theta + \omega +\Omega -\lambda_{31})\\
 & \qquad + \Bigl(\frac{15}{16} \sin^2 i (1-3\cos i)-\frac{3}{4} (1-\cos i) \Bigr) \frac{11 e^2}{8} \cos (M-\theta - \omega +\Omega -\lambda_{31})\Bigr\}\\
& \ + \frac{\mu_E R_E^3 J_{32}}{a^4} \Bigl\{ -\frac{15}{8} \sin i (1+\cos i)^2 e \sin (2M-2\theta +3 \omega +2\Omega -2\lambda_{32}) \\
& \qquad + \frac{45}{8} \sin i (1-2\cos i-3 \cos^2 i) e \sin (2M-2\theta + \omega +2\Omega -2\lambda_{32})\Bigr\}\\
& \ + \frac{\mu_E R_E^3 J_{33}}{a^4} \Bigl\{ \frac{15}{8} (1+\cos i)^3 (1-6 e^2) \cos [3(M-\theta + \omega +\Omega -\lambda_{33})] \\
& \qquad + \frac{45}{8} \sin^2 i (1+\cos i) \frac{53 e^2}{8} \cos (3M-3\theta + \omega+ 3\Omega -3\lambda_{33})\Bigr\}\\
& \ + \frac{\mu_E R_E^4 J_{41}}{a^5} \Bigl\{ \Bigl(\frac{35}{16} \sin^3 i (1+2\cos i)-\frac{15}{8}(1+\cos i)\sin i \Bigr) \frac{e}{2} \sin (M-\theta +2 \omega +\Omega -\lambda_{41}) \\
& \qquad + \cos i \Bigl(\frac{15}{4} \sin i -\frac{105}{16} \sin^3 i \Bigr) \frac{5e}{2} \sin (M-\theta + \Omega -\lambda_{41})\Bigr\}\\
\end{split}
\end{equation}

\begin{equation}\label{term11}
\begin{split}
& \ + \frac{\mu_E R_E^4 J_{42}}{a^5} \Bigl\{ -\frac{105}{32} \sin^2 i (1+\cos i)^2\ \frac{e^2}{2} \cos[2 (M-\theta +2 \omega +\Omega -\lambda_{42})] \\
&  \qquad +  \Bigl(\frac{105}{8} \sin^2 i \cos i (1+\cos i) -\frac{15}{8} (1+\cos i)^2 \Bigr) (1+e^2) \cos[2 (M-\theta + \omega+ \Omega -\lambda_{42})]\\
& \qquad + \Bigl(\frac{105}{16} \sin^2 i (1-3\cos^2 i) -\frac{15}{4} \sin^2 i \Bigr) 5e^2 \cos[2 (M-\theta + \Omega -\lambda_{42})]\Bigr\}\\
& \ + \frac{\mu_E R_E^4 J_{43}}{a^5} \Bigl\{ \frac{105}{16} \sin i (1+\cos i)^3 \Bigl(-\frac{3e}{2} \Bigr) \sin (3M-3\theta +4 \omega +3\Omega -3\lambda_{43}) \\
& \qquad + \frac{105}{8} \sin i (1-3\cos^2 i-2 \cos^3 i) \Bigl(\frac{9e}{2} \Bigr) \sin (3M-3\theta + 2\omega +3\Omega -3\lambda_{43})\Bigr\}\\
& \ + \frac{\mu_E R_E^4 J_{44}}{a^5} \Bigl\{ \frac{105}{16} (1+\cos i)^4 \Bigl(1-11 e^2 \Bigr) \cos [4(M-\theta + \omega +\Omega -\lambda_{44})] \\
& \qquad + \frac{105}{4} \sin^2 i (1+ \cos i)^2 \Bigl(\frac{53e^2}{4} \Bigr) \cos (4M-4\theta + 2\omega +4\Omega -4\lambda_{44})\Bigr\}\ .
\end{split}
\end{equation}

From Table~\ref{table:J} we see that the leading coefficients are $J_{22}$, $J_{31}$. Therefore, with reference to
\equ{term11}, let us introduce the following notation:
\begin{equation}
\begin{split}\label{t11_terms}
& \mathfrak{T}_1=g_1 (L,G,H) \cos [2 (\lambda -\lambda_{22})]\\
& \mathfrak{T}_2=g_2 (L,G,H) \cos [2 (\lambda -\omega -\lambda_{22})]\\
& \mathfrak{T}_3=g_3 (L,G,H) \cos (\lambda -\lambda_{31})\ ,
\end{split}
\end{equation}
where
\begin{equation}
\begin{split}\label{g11_terms}
& g_1= \frac{\mu_E R_E^2 J_{22}}{a^3} \frac{3}{4} (1+\cos i)^2 (1-\frac{5}{2}e^2)\\
& g_2= \frac{\mu_E R_E^2 J_{22}}{a^3} \frac{27}{8} e^2 \sin^2 i\\
& g_3= \frac{\mu_E R_E^3 J_{31}}{a^4} \Bigl\{ \frac{15}{16}
\sin^2 i (1+3\cos i) -\frac{3}{4}(1+\cos i) \Bigr\} (1+2 e^2)\ ,
\end{split}
\end{equation}
and
$$\lambda=M-\theta+\omega+\Omega,$$ is the so--called {\it stroboscopic mean node}.

The magnitude of each term in \equ{term11} varies with the eccentricity and the inclination. In order to compare the effects produced by the terms of $R_{earth}^{res1:1}$ and to provide an analytical argument for the numerical results which will be presented in Section~\ref{sec:cartography}, we introduce the following heuristic definition of \sl dominant \rm term.

\begin{definition} \label{def_dominant}
For given values of $(a,e,i)$, equivalently for given values of $(L,G,H)$,
we say that a specific term, say $\mathfrak{T}_k$ for some $k\in\mathbb{Z}$, of the expansion of
$R_{earth}^{res1:1}$ is {\it dominant}  with respect to the other harmonic terms of the resonant part, if the magnitude of $|g_k(L,G,H)|$ is greater than the magnitude of any other term of the expansion.
\end{definition}


\begin{table}[h]
\begin{tabular}{|c|c|c|c|c|c|}
  \hline
  $n$ & $m$ & $C_{nm}$ & $S_{nm}$ & $J_{nm}$ & $\lambda_{nm}$ \\
  \hline
  2 & 0 & -1082.6261& 0& 1082.6261& 0 \\
  2 & 1 & -0.000267& 0.0017873& 0.001807& $-81_{\cdot}^{\circ}5116$ \\
2 & 2 & 1.57462& -0.90387& 1.81559& $75_{\cdot}^{\circ}0715$ \\
 3 & 0 & 2.53241& 0& -2.53241& 0 \\
3 & 1 & 2.19315& 0.268087& 2.20947& $186_{\cdot}^{\circ}9692$ \\
3 & 2 & 0.30904& -0.211431& 0.37445& $72_{\cdot}^{\circ}8111$ \\
3 & 3 & 0.100583& 0.197222& 0.22139& $80_{\cdot}^{\circ}9928$ \\
4 & 0 & 1.6199& 0& -1.619331& 0 \\
 4 & 1 & -0.50864& -0.449265& 0.67864& $41_{\cdot}^{\circ}4529$ \\
4 & 2 & 0.078374& 0.148135& 0.16759& $121_{\cdot}^{\circ}0589$ \\
 4 & 3 & 0.059215& -0.012009& 0.060421& $56_{\cdot}^{\circ}1784$ \\
 4 & 4 & -0.003983& 0.006525& 0.007644& $-14_{\cdot}^{\circ}6491$ \\
  \hline
 \end{tabular}
 \vskip.1in
 \caption{Values of $C_{nm}$, $S_{nm}$, $J_{nm}$ (in units of $10^{-6}$), computed from \cite{EGM2008}.}\label{table:J}
\end{table}

\begin{figure}[h]
\centering
\vglue-1.1cm
\includegraphics[width=6truecm,height=5truecm]{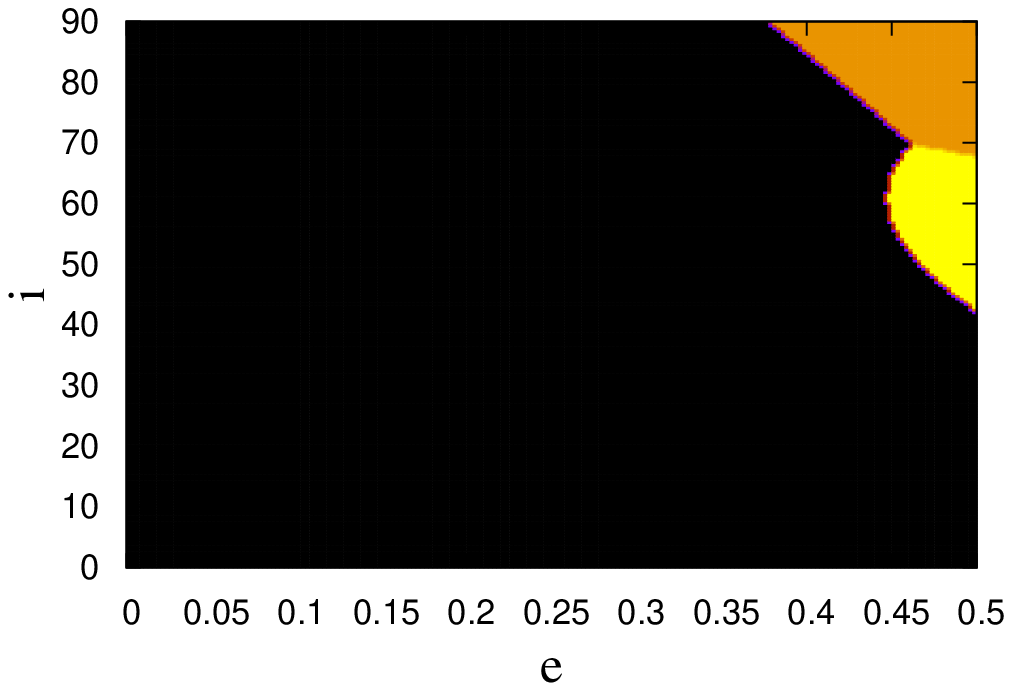}
\includegraphics[width=6truecm,height=5truecm]{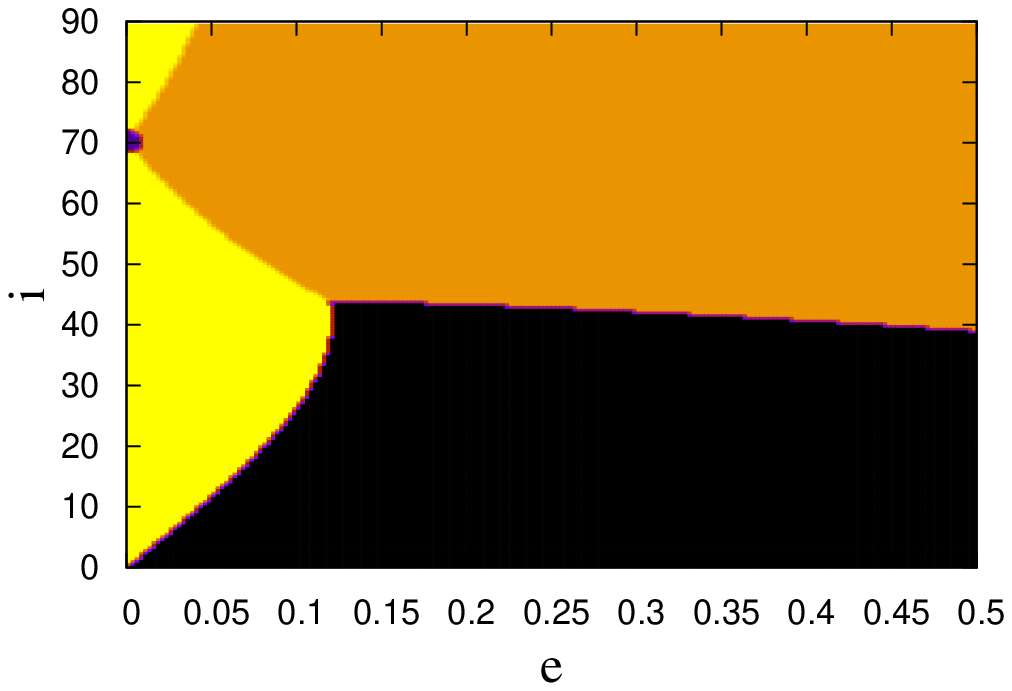}
\vglue0.3cm
\caption{Left: dominant terms in $R_{earth}^{res1:1}$ as a function of
$(e,i)$: $\mathfrak{T}_1$--black, $\mathfrak{T}_2$--brown,
$\mathfrak{T}_3$--yellow, where $\mathfrak{T}_1$, $\mathfrak{T}_2$, $\mathfrak{T}_3$ are defined in \equ{t11_terms}.
Right: dominant terms in $R_{earth}^{res2:1}$ as a function of
$(e,i)$: $t_1$--black, $t_2$--brown,
$t_3$--yellow, where $t_1$, $t_2$, $t_3$ are defined in \equ{t1t2t3}.
The blue dot at $(e,i)\simeq(0^o,70^o)$ corresponds to a fourth degree harmonic term.} \label{big11}
\end{figure}

We proceed to provide an analysis of the dominant terms, which turns out to be a simple, but essential tool to perform fast computations.
In particular, we underline which are the most important terms to be taken into account in the expansion of the geopotential, thus avoiding the integration of very long expressions, like those given in \equ{term11}.
With reference to \equ{t11_terms}, in Figure~\ref{big11} left we represent the index of the dominant term
as a function of eccentricity and inclination, where the colors are
set as follows: black means that $\max \{|g_1|, |g_2|, |g_3|\}=|g_1|$,
brown shows that $|g_2|$ has the highest value, yellow expresses the
fact that $\mathfrak{T}_3$ dominates. From the analysis of Figure~\ref{big11} left we conclude that $
\mathfrak{T}_1$ dominates in almost all regions of the plane, except for some high eccentricities and
inclinations.
A refined analysis shows that for small inclinations and eccentricities the magnitude of  $\mathfrak{T}_1$ is much greater than the magnitude of any other term in the expansion. This will be the main reason for getting a pendulum like behavior with the stable point located at $\lambda=\lambda_{22}$, see Remark~\ref{rem:equil} and compare with Figure~\ref{fligeo} below. In conclusion, for given values
of eccentricity and inclination, Figure~\ref{big11} provides the dominant term
of the expansion \equ{term11}. Of course the analysis can be extended by considering more terms of the expansion \equ{term11},
but we limit our discussion to the most relevant ones leaded by the coefficients $J_{22}$, $J_{31}$, whose size is the
highest one, as it is shown in Table~\ref{table:J}.

\subsubsection{Resonance 2:1.}  The 2:1 resonant terms correspond to
$2(n-2p+q)=m$, $m>0$. In this case we extend the computations up to the fourth
order in the eccentricity, since in the GPS region there are satellites
(e.g., Molniya satellites) with very high values of the eccentricity, say
$e\in(0.5,0.75)$. Thus, we obtain the following expansion:
\begin{equation}
\begin{split}
& R_{earth}^{res2:1}\cong \frac{\mu_E R_E^2 J_{22}}{a^3} \Bigl\{\frac{3}{4} (1+\cos i)^2 \Bigl( -\frac{e}{2}+\frac{e^3}{16}\Bigr) \cos (M-2\theta +2\omega +2\Omega -2\lambda_{22})\nonumber\\
& \qquad + \frac{3}{2} \sin^2 i \Bigl(\frac{3}{2}e +\frac{27}{16}
e^3\Bigr)  \cos (M-2\theta +2\Omega -2\lambda_{22}) \nonumber\\
& \qquad + \frac{3}{4} (1-\cos i)^2 \frac{67 e^3}{48}  \cos
(M-2\theta -2\omega +2\Omega -2\lambda_{22})
\Bigr\}\nonumber\\
& \ + \frac{\mu_E R_E^3 J_{32}}{a^4} \Bigl\{ \frac{15}{8} \sin i (1+\cos i)^2 \Bigl( \frac{e^2}{8}+\frac{e^4}{48} \Bigr) \sin (M-2\theta +3 \omega +2\Omega -2\lambda_{32}) \nonumber\\
& \qquad + \frac{15}{8} \sin i (1-2\cos i-3 \cos^2 i) \Bigl(1+2e^2+\frac{239 e^4}{64}\Bigr) \sin (M-2\theta + \omega +2\Omega -2\lambda_{32}) \nonumber\\
&  \qquad - \frac{15}{8} \sin i (1+2\cos i-3 \cos^2 i) \Bigl(\frac{ 11 e^2}{8} +\frac{49 e^4}{16} \Bigr) \sin (M-2\theta - \omega +2\Omega -2\lambda_{32}) \nonumber\\
& \qquad - \frac{15}{8} \sin i (1- \cos i)^2 \frac{131 e^4}{128}
\sin (M-2\theta - 3\omega +2\Omega -2\lambda_{32})
 \Bigr\}\nonumber\\
 & \ + \frac{\mu_E R_E^4 J_{42}}{a^5} \Bigl\{ \frac{105}{32} \sin^2
i  (1+\cos i)^2 \frac{19 e^3}{48}
 \cos(M-2\theta + 4\omega+ 2\Omega -2\lambda_{42})\nonumber\\
  & \qquad +\Bigl(\frac{105}{8} \sin^2 i \cos i (1+\cos i) -\frac{15}{8} (1+\cos i)^2 \Bigr) \Bigl(\frac{e}{2}+\frac{33 e^3}{16} \Bigr) \cos(M-2\theta + 2\omega+ 2\Omega -2\lambda_{42})\nonumber\\
& \qquad + \Bigl(\frac{105}{16} \sin^2 i (1-3\cos^2 i) -\frac{15}{4} \sin^2 i \Bigr) \Bigl( \frac{5e}{2}+\frac{135 e^3}{16} \Bigr) \cos(M-2\theta + 2\Omega -2\lambda_{42})\nonumber\\
& \qquad - \Bigl(\frac{105}{8} \sin^2 i \cos i (1-\cos i)+\frac{15}{8} (1- \cos i)^2 \Bigr) \frac{49 e^3}{48}
\cos(M-2\theta -2 \omega + 2\Omega -2\lambda_{42})\Bigr\}\nonumber
\end{split}
\end{equation}

\begin{equation}
\begin{split}
& \ + \frac{\mu_E R_E^4 J_{44}}{a^5} \Bigl\{ \frac{105}{16} (1+\cos i)^4 \Big(\frac{e^2}{2}-\frac{e^4}{3} \Bigr) \cos [2(M-2\theta + 2\omega +2\Omega -2\lambda_{44})] \nonumber\\
& \qquad + \frac{105}{4} \sin^2 i (1+ \cos i)^2 \Bigl(1+e^2 +\frac{65 e^4}{16} \Bigr) \cos [2(M-2\theta + \omega +2\Omega -2\lambda_{44})]\nonumber\\
& \qquad + \frac{315}{8} \sin^4 i \Bigl( 5 e^2+\frac{155 e^4}{12} \Bigr) \cos [2(M-2\theta + 2\Omega -2\lambda_{44})]\nonumber\\
& \qquad + \frac{105}{4} \sin^2 i (1-\cos i)^2 \frac{67 e^4}{48} \cos [2(M-2\theta -\omega + 2\Omega -2\lambda_{44})]\Bigr\}\ .\nonumber
 \end{split}
\end{equation}

\vskip.1in

Definition~\ref{def_dominant} of \sl dominant terms \rm can be straightforwardly extended to any resonance. In particular,
by comparing the magnitude of each
harmonic in $R_{earth}^{res2:1}$, we note that the
most important terms in $R_{earth}^{res2:1}$ are $t_1$, $t_2$ and $t_3$, where
\begin{equation}
\begin{split}\label{t1t2t3}
& t_1= \frac{\mu_E R_E^2 J_{22}}{a^3} \Bigl\{\frac{3}{4} (1+\cos i)^2 \Bigl( -\frac{e}{2}+\frac{e^3}{16}\Bigr) \cos [2(\lambda +\frac{\omega}{2} -\lambda_{22})]\Bigr\}\\
& t_2=\frac{\mu_E R_E^2 J_{22}}{a^3} \Bigl\{\frac{3}{2} \sin^2 i
\Bigl(\frac{3}{2}e+\frac{27}{16}e^3 \Bigr)  \cos [2(\lambda -\frac{\omega}{2}-\lambda_{22})] \Bigr\} \\
& t_3= \frac{\mu_E R_E^3 J_{32}}{a^4} \Bigl\{ \frac{15}{8} \sin i
(1-2\cos i-3 \cos^2 i) \Bigl(1+2e^2+\frac{239}{64}e^4 \Bigr) \sin[2 (\lambda -\lambda_{32})]
\Bigr\}\\
\end{split}
\end{equation}
and $\lambda$ is the stroboscopic mean node defined in the case of the $2:1$ resonance by
\begin{equation}\label{lambda_definition_2:1}
\lambda=\frac{1}{2}\Bigl[M+\omega-2(\theta-\Omega)\Bigr].
\end{equation}

The right panel of Figure~\ref{big11} is the analogue of the left one for the 2:1 resonance and it shows the dominant term
as a function of eccentricity and inclination, where the colors are
set as follows: black means that $t_1$ dominates,
brown shows that $t_2$ has the highest magnitude, yellow expresses the
fact that $t_3$ is the  dominant term.

In contrast to the $1:1$ resonance, where a single term, precisely $\mathfrak{T}_1$, prevails over a large region of the domain, we have a different balance between the three terms. In particular, for
many eccentricities and inclinations the magnitude of $\mathfrak{T}_1$ in Figure~\ref{big11} left was greater than the magnitude of any other term of $R^{res1:1}_{earth}$, while for the $2:1$ resonance we do not have such a prominent term. Each of the terms $t_1$, $t_2$ or $t_3$ dominates (in the sense of Definition~\ref{def_dominant}) in some specific region of the domain $(e,i)$. However, for the 2:1 resonance if one term dominates, it does not mean that its magnitude is \sl much \rm greater than the magnitude of the other two terms. Indeed, there are large regions of the domain, where $t_1$, $t_2$ and $t_3$ are comparable in magnitude. In this respect the results for the 2:1 resonance are much different than in the case of the $1:1$ resonance.

There is a special case where a fourth degree harmonic term dominates. This case stems from the fact that all
terms (except $t_3$ and a $J_{44}$--term) are of order $O(e)$ and therefore they are zero for $e=0$; however, since $t_3$ is zero for $i=70.53^o$,
a blue dot in Figure~\ref{big11} right, having the coordinates $(e=0, i=70.53^o)$,
indicates that the fourth  degree harmonic term is dominant.

However, this is a singular case and usually one of the terms $t_1$, $t_2$ or $t_3$ dominates. Excluding this singular case, it is important to point out that for small and moderate eccentricities, $t_1$, $t_2$, $t_3$ are dominant not just in the sense of the above definition, but at least one of them has its magnitude much larger than the magnitude of any other term of the expansion. For such eccentricities, all major features of the resonance patterns are due to
these three terms and any other harmonic plays a less relevant role.

\begin{remark}\label{rem:equil}
We add an elementary remark about the determination of the
equilibrium points associated to the resonances and their
stability. In the case of the 1:1 resonance, with reference to
\equ{t11_terms} and \equ{g11_terms}, we write the equations of
motion associated, e.g. to the term $g_1$, for the variables
$\lambda$ and $L$ as
\beqa{lambdaL}
\dot\lambda&=&\frac{\mu^2}{L^3}-\dot\theta+\frac{\partial g_1}{\partial L}\cos(2\lambda-2\lambda_{22})\nonumber\\
\dot L&=&2 g_1 \sin(2\lambda-2\lambda_{22})\ .
\eeqa
The equilibrium points are given by $\lambda=\lambda_{22}$ or $\lambda=\lambda_{22}+\pi/2$ (so that $\dot L=0$) for a suitable $L=L_0$
(so that $\dot\lambda=0$).
The stability is computed by looking at
the eigenvalues of the derivative of the vector field \equ{lambdaL} at the equilibrium points.
For $\lambda=\lambda_{22}$ we obtain that the eigenvalues
are purely imaginary and therefore the equilibrium is stable; on the contrary, for
$\lambda=\lambda_{22}+\pi/2$ we obtain real eigenvalues, thus leading to a linear instability of the equilibrium position.
A similar discussion holds for the other terms of the 1:1 resonance as well as for
the determination of the equilibria associated to other resonances.
\end{remark}

\section{A measure of the amplitude of resonant islands}\label{sec:amplitude}

In this section we provide an elementary computation to estimate the amplitude arising around a given $p:q$ resonance.
We stress that the formulae we shall derive can be implemented without any computational effort and yet provide a very reliable estimate
of the size of the resonant islands. This analysis, complemented with that of the dominant terms presented in Section~\ref{sec:ham}, provides a
fast description of the dynamics in the neighborhood of a resonance.

Having denoted by $R^{sec}_{earth}$ the secular part of the geopotential (see \equ{Rsec}), we consider the Hamiltonian corresponding
to the $p:q$ resonance (for some $p$, $q$ coprime integers), where the non--resonant terms have been averaged out.
We expand the resonant part in Fourier series up to finite orders, say $N_1$, $N_2$, $N_3$, and we denote by
$R_{\underline k}^{(p,q)}$ the Fourier coefficients of the resonant part. The resulting resonant Hamiltonian is thus given by
\beqa{Hres}
&&\mathcal{H}_{res}^{p:q}(L,G,H,qM-p\theta,\omega,\Omega)=-{{\mu_E^2}\over {2L^2}}+R^{sec}_{earth}(L,G,H,\omega)\nonumber\\
&+&\sum_{k_1=1}^{N_1}\sum_{k_2=1}^{N_2}\sum_{k_3=1}^{N_3}
R_{\underline k}^{(p,q)}(L,G,H)\ cs(k_1(q M-p\theta)+k_2\omega+k_3\Omega)\ ,
\eeqa
where $cs$ could be either cosine or sine. We recall that the $p:q$ resonance corresponds to the relation
$q\dot M-p\dot\theta=0$;
taking into account that $\dot\theta=1$ (in normalized units) and denoting by $L_{res}$ the resonant value of the Delaunay action $L$, we have
$$
\dot M={p\over q}={{\mu_E^2}\over {L_{res}^3}}\ ,
$$
so that the resonant value $L_{res}$ is obtained as
\beq{H5}
L_{res}=\Big({q\over p}\ \mu_E^2\Big)^{1\over 3}\ .
\eeq
Using that $a=L^2/\mu_E$, we get that the resonant value of the semimajor axis is
$$
a_{res}=\Big({q\over p}\Big)^{2\over 3}\ \mu_E^{1\over 3}\ .
$$

Let us proceed to expand \equ{Hres} around $L_{res}$ up to second
order: \beqa{H1}
&&\mathcal{H}_{res}^{p:q}(L,G,H,qM-p\theta,\omega,\Omega)=-{{\mu_E^2}\over
{2L_{res}^2}}+{{\mu_E^2}\over {L_{res}^3}}\,(L-L_{res})
-{{3\mu_E^2}\over {2L_{res}^4}}\,(L-L_{res})^2\nonumber\\
&+&R^{sec}_{earth}(L_{res},G,H,\omega)+R^{sec,L}_{earth}(L_{res},G,H,\omega)\,(L-L_{res})\nonumber\\
&+&{1\over 2}R^{sec,LL}_{earth}(L_{res},G,H,\omega)\,(L-L_{res})^2\nonumber\\
&+&\sum_{k_1=1}^{N_1}\sum_{k_2=1}^{N_2}\sum_{k_3=1}^{N_3}
R_{\underline k}^{(p,q)}(L_{res},G,H)\ cs(k_1(q M-p\theta)+k_2\omega+k_3\Omega)\ ,
\eeqa
where $R^{sec,L}_{earth}$, $R^{sec,LL}_{earth}$ denote first and second derivatives with respect to $L$ of $R^{sec}_{earth}$.
Setting $\Lambda\equiv L-L_{res}$ and neglecting constant terms as well as the term $R^{sec}_{earth}(L_{res},G,H,\omega)$, we rewrite
\equ{H1} retaining only the largest term in the resonant Hamiltonian:
\beq{H2}
\mathcal{H}_{res}^{p:q,max}(\Lambda,G,H,qM-p\theta,\omega,\Omega)=\alpha\Lambda-\beta\Lambda^2+A\ cs(k_1^{max}(qM-p\theta)
+k_2^{max}\omega+k_3^{max}\Omega)\ ,
\eeq
where we have introduced the quantities
\beqa{H4}
\alpha&=&\alpha(L_{res},G,H,\omega)\equiv{{\mu_E^2}\over {L_{res}^3}}+R^{sec,L}_{earth}(L_{res},G,H,\omega)\nonumber\\
\beta&=&\beta(L_{res},G,H,\omega)\equiv{{3\mu_E^2}\over {2L_{res}^4}}-{1\over 2}R^{sec,LL}_{earth}(L_{res},G,H,\omega)\nonumber\\
A&=&A(L_{res},G,H)\equiv\max_{k_1,k_2,k_3} R_{\underline k}^{(p,q)}(L_{res},G,H)
\eeqa
and ${\underline k}^{max}\equiv(k_1^{max},k_2^{max},k_3^{max})$ is the index at which the maximum defining $A$ is reached.

Since the second term in $\beta$ is much smaller than the first term, in practical computations we replace $\beta$ with
\beq{H4bis}
\beta\equiv{{3\mu_E^2}\over {2L_{res}^4}}\ .
\eeq
Let us find the quantities $B$, $C$ such that
$$
\alpha\Lambda-\beta\Lambda^2=-(B+C\Lambda)^2+B^2\ ;
$$
we immediately find
$$
B=-{\alpha\over {2\sqrt{\beta}}}\ ,\qquad C=\sqrt{\beta}\ .
$$
Neglecting again constant terms, we write \equ{H2} as
\beq{H3}
\mathcal{H}_{res}^{p:q}(\Lambda,G,H,qM-p\theta,\omega,\Omega)=-(B+C\Lambda)^2+A\ cs(k_1^{max}(qM-p\theta)+k_2^{max}\omega+k_3^{max}\Omega)\ .
\eeq
Next, we define $\Gamma\equiv B+C\Lambda$, so that \equ{H3} becomes
$$
\mathcal{H}_{res}^{p:q}(\Gamma,G,H,qM-p\theta,\omega,\Omega)=-{({\sqrt{2}\Gamma)^2}\over 2}+A\ cs(k_1^{max}(qM-p\theta)+k_2^{max}\omega+k_3^{max}\Omega)\ .
$$
Then, we obtain that the excursion in $\Gamma$, say $\Delta\Gamma$, is equal to (compare with \cite{LL})
$\Delta\Gamma=\sqrt{2A}$, which provides
$$
\Delta\Lambda={1\over C}\sqrt{2A}=\sqrt{{2A}\over \beta}\ .
$$
Going back to the Delaunay action $L$, we get $\Delta L=\sqrt{2A/\beta}$;
taking into account that $a=L^2/\mu_E$, we have
$$
\Delta a={L^2\over \mu_E}-{L_{res}^2\over \mu_E}={1\over \mu_E}\ (\Delta L^2+2L_{res}\ \Delta L)\ .
$$
Therefore, we obtain that the full amplitude of the $p:q$ resonant island, measured in terms of the semimajor axis, is given by
$$
2\ \Delta a={2\over \mu_E}\ (\Delta L^2+2L_{res}\ \Delta L)={2\over \mu_E}\ ({{2A}\over \beta}+2L_{res}\ \sqrt{{2A}\over \beta})
$$
with $A$, $\beta$ as in \equ{H4}, \equ{H4bis} and $L_{res}$ as in \equ{H5}.

\begin{figure}[h]
\centering
\vglue-0.9cm
\includegraphics[width=6truecm,height=5truecm]{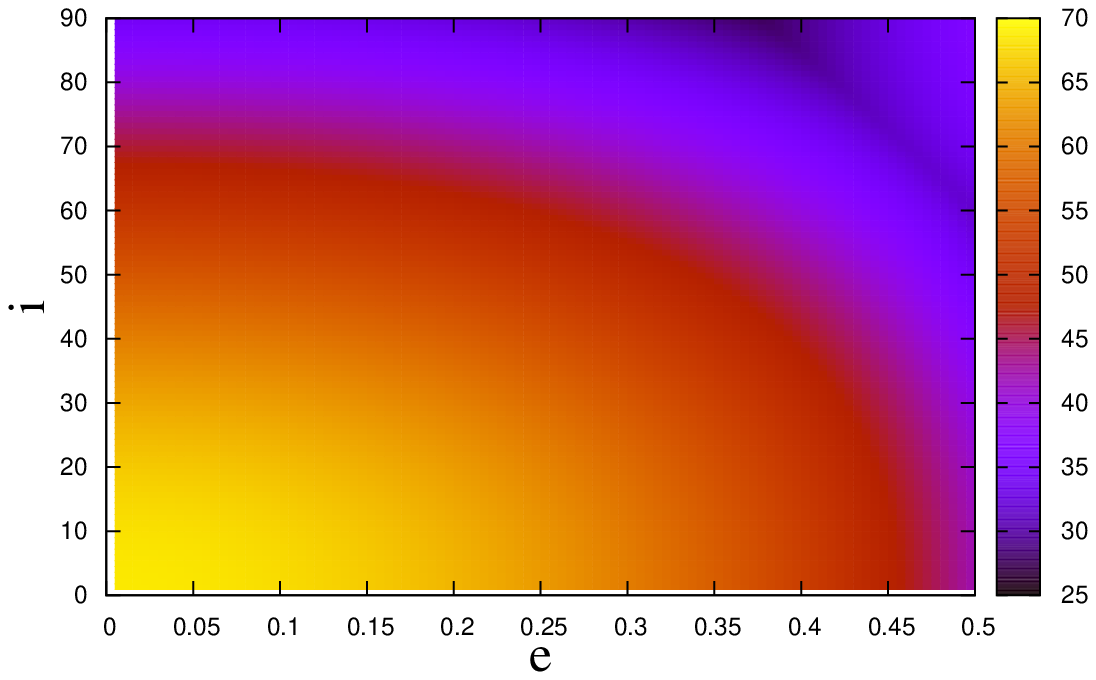}
\includegraphics[width=6truecm,height=5truecm]{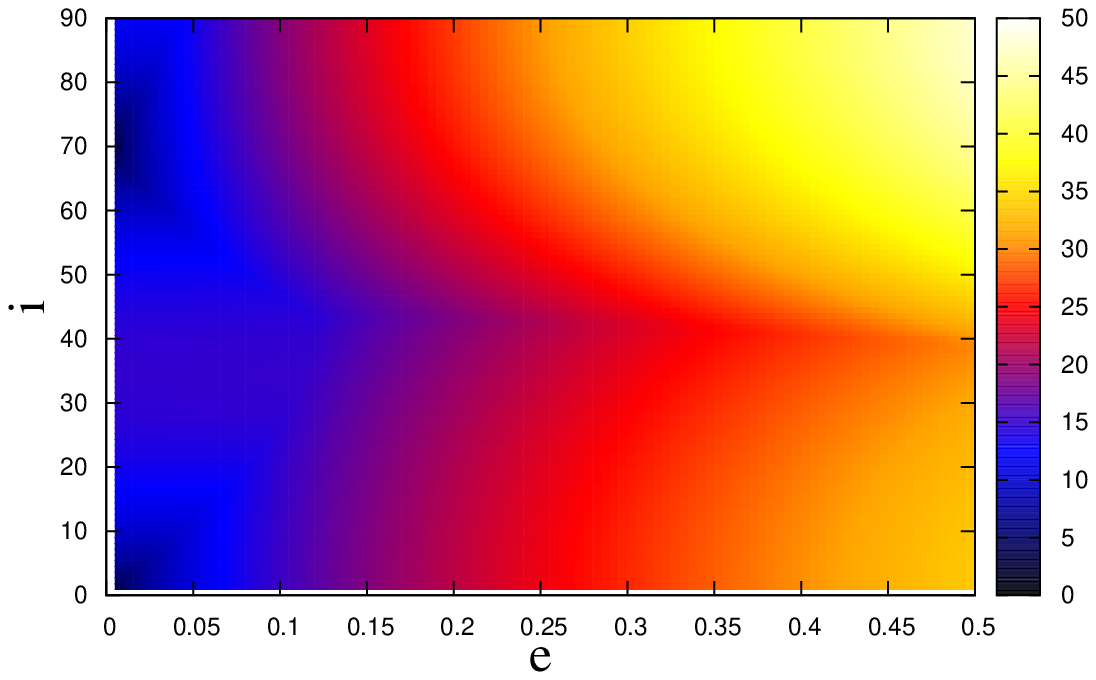}
\vglue0.3cm
\caption{The amplitude of the resonances for different values of the eccentricity
and the inclination; the color bar provides the measure of the amplitude in kilometers.
Left: the 1:1 resonance; right: the 2:1 resonance. } \label{fig:amplitude}
\end{figure}

Figure~\ref{fig:amplitude} provides the amplitudes of the 1:1 and 2:1 resonances as the eccentricity varies between
0 and 0.5, while the inclination ranges between $0^o$ and $90^o$ (we took $\omega=0$ and $\Omega=0$). The color bar indicates
the size of the amplitude in kilometers. All values are in agreement with the amplitudes shown in the figures of
the forthcoming Section~\ref{sec:cartography}, either when the results are obtained through
Hamiltonian or Cartesian formalism. The method for computing the resonant amplitudes
extends very well to other resonances of different order, even if the size is small (see \cite{CGminor}).
Notice that the peculiar behavior in Figure~\ref{fig:amplitude}, right panel, at about $i=40^o$ is due to the fact that
such inclination corresponds to the boundary of the regions where the terms $t_1$ and $t_2$ are dominant.
The analytical determination of this value will be given through equation \equ{i0} in Section~\ref{sec:cart21}.

\section{Cartography}\label{sec:cartography}
In this section we describe the results about the analysis of the
1:1 and 2:1 resonances using the FLIs, whose definition is
provided below in Section~\ref{sec:FLI}. Each region shows a
peculiar cartography, where the main zones of resonant, stable and
chaotic motions are highlighted.

\subsection{Fast Lyapunov Indicator}\label{sec:FLI}
The FLI can be defined as the largest Lyapunov characteristic exponent
at a fixed time, say $t=T$ (\cite{froeslega}). More precisely, let us consider the $n$--dimensional differential system
$$
\dot{\bf{x}}=\bf{f}(\bf{x})\ ,
$$
where ${\bf{x}}\in {\mathbb{R}}^{n}$ and the $n$--dimensional vector function ${\bf{f}}:{\mathbb{R}}^{n}\rightarrow {\mathbb{R}}^{n}$ represents the vector field. Let the corresponding variational equations be written as
$$
\dot {\bf{v}}=\Big({{\partial \bf{f}(\bf{x})} \over {\partial \bf{x}}}\Big)\ \bf{v}\ ,
$$
where $\bf{v}$ is an $n$-dimensional vector.

Given the initial conditions ${\bf{x}}(0) \in \mathbb{R}^{n}$,
${\bf{v}}(0) \in \mathbb{R}^{n}$, the FLI at the time $T>0$ is given by
$$
{\rm FLI}({\bf{x}}(0), {\bf{v}}(0), T) \equiv \sup _{0 < t\leq T} \log ||{\bf{v}}(t)||\  .
$$
In the following sections (Section~\ref{sec:cart11} for the 1:1 resonance and Section~\ref{sec:cart21} for the
2:1 resonance) we present results in the plane of coordinates or in the parameter plane, providing the value of the FLI through a color scale, where darker colors will denote a regular dynamics, either periodic or quasi--periodic, while lighter colors will denote chaotic motions. We remark that in each plot the color scale may be different.

The results of Sections~\ref{sec:cart11} and \ref{sec:cart21} will provide
information on the regular or chaotic character of the dynamics, on the dependence of the resonances on the parameters,
on the location of the equilibrium points. Our study will be mainly based on the Hamiltonian formulation for the 1:1 and 2:1
resonances, but we will also provide results using the Cartesian approach including, beside the geopotential, the effects
of Sun, Moon and solar radiation pressure.

\subsection{Cartography of the 1:1 resonance}\label{sec:cart11}

\begin{figure}
\centering
\vglue-0.9cm
\includegraphics[width=6truecm,height=5truecm]{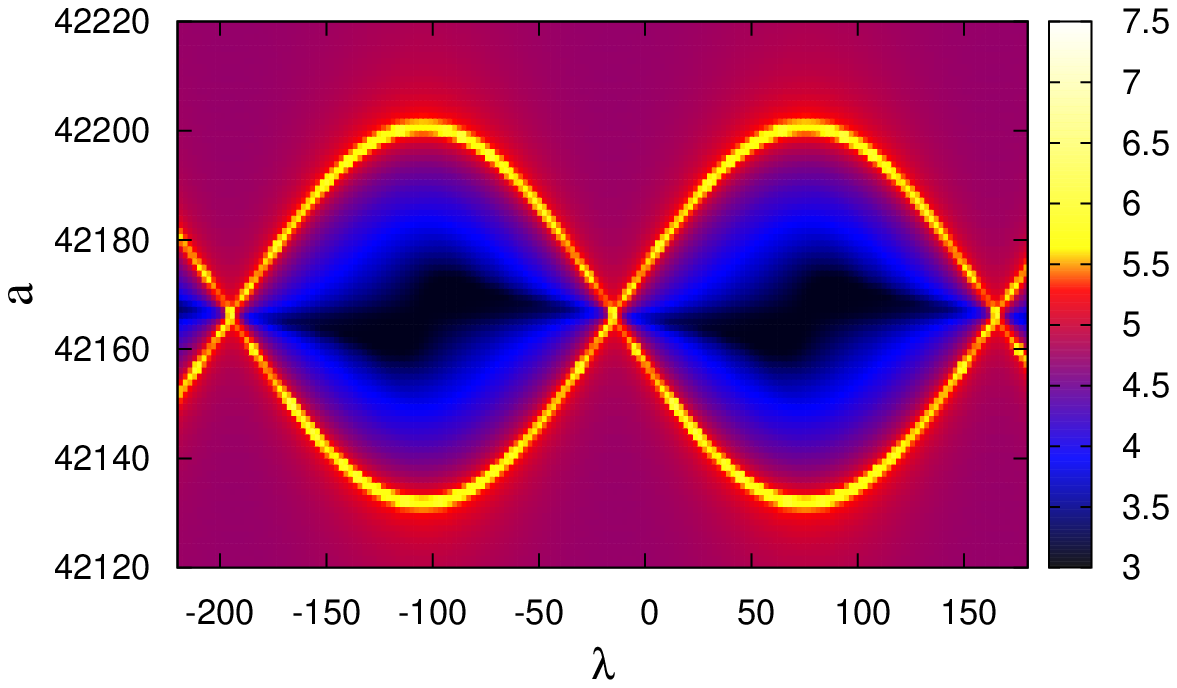}
\includegraphics[width=6truecm,height=5truecm]{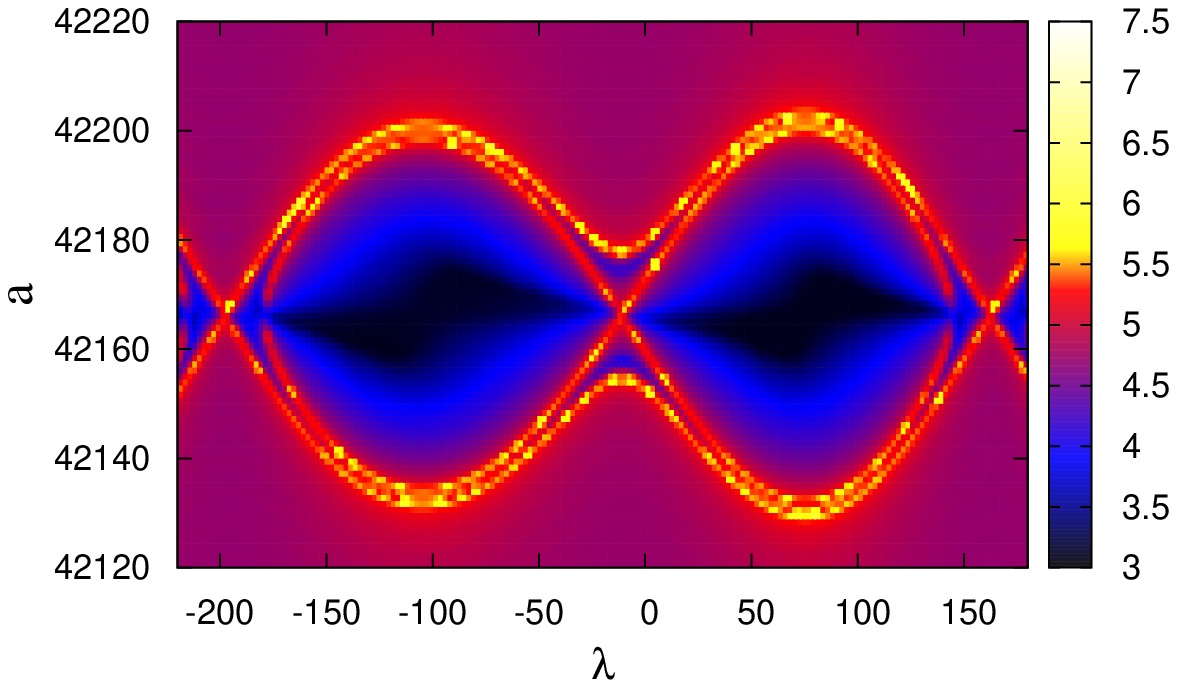}\\
\vglue-0.8cm
\includegraphics[width=6truecm,height=5truecm]{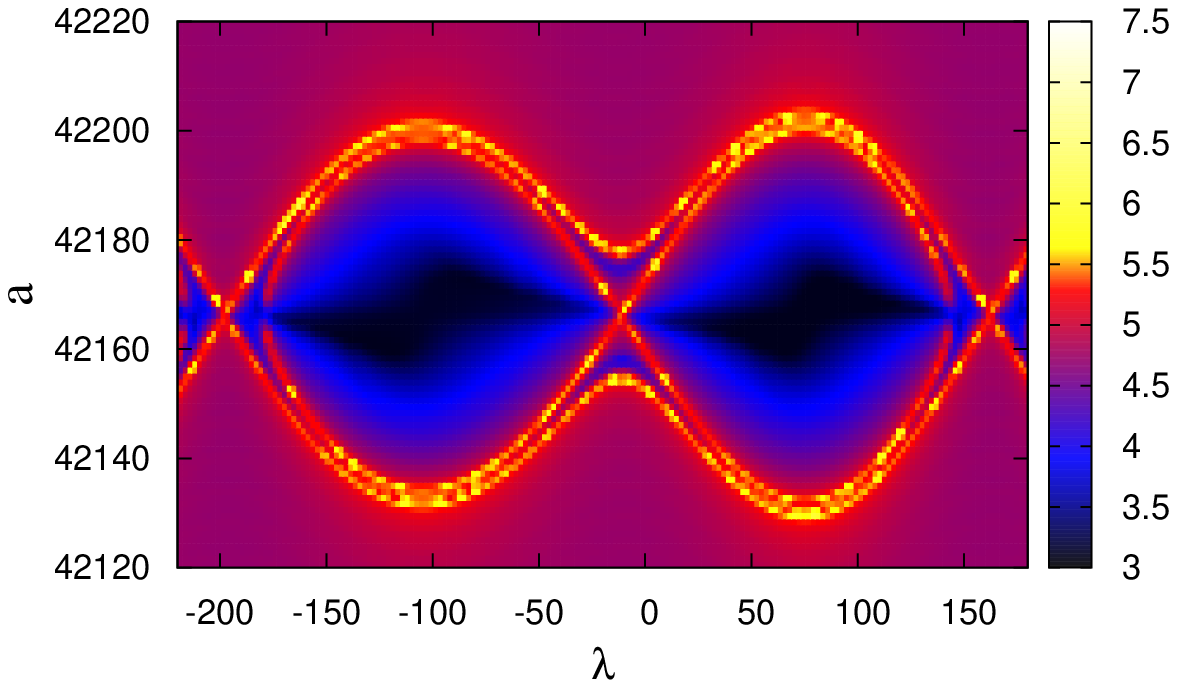}
\includegraphics[width=6truecm,height=5truecm]{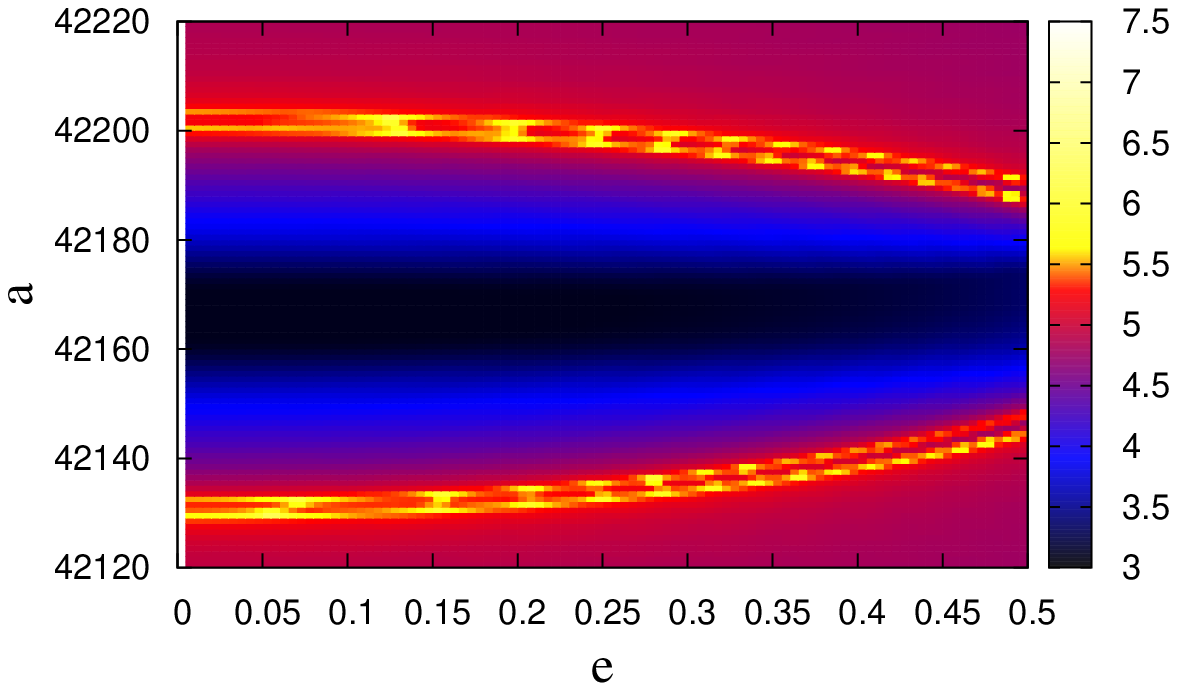}
\vglue0.3cm
\caption{FLI (using Hamilton's equations) for the GEO 1:1 resonance for $e=0.005$, $i=0^o$, $\omega=0^o$, $\Omega=0^o$ under the effects of
the $J_2$ and $J_{22}$ terms (top left); all harmonics up to degree and order $n=m=3$ (top right);
all harmonics up to $n=m=4$ (bottom left). The bottom right panel yields the FLI for $i=0^o$, $\lambda=75.07^o$
in the $(e,a)$ plane under the effects of all harmonics up to $n=m=4$. }
\label{fligeo}
\end{figure}

\begin{figure}
\centering
\vglue-0.9cm
\includegraphics[width=6truecm,height=5truecm]{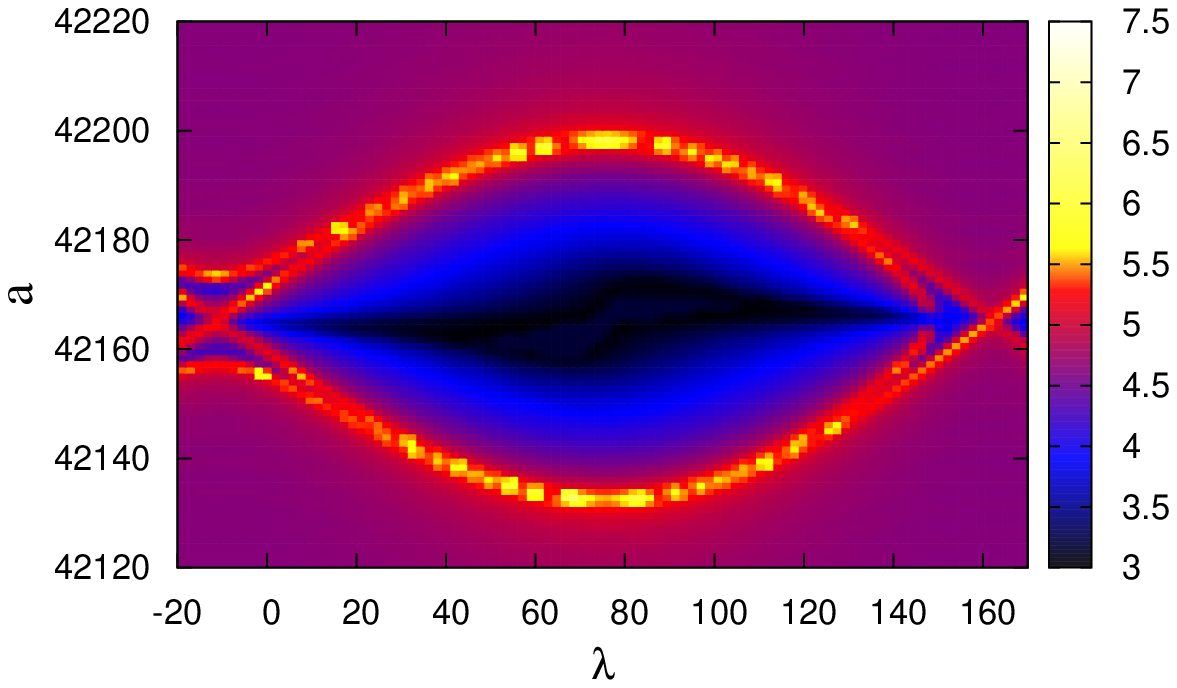}
\includegraphics[width=6truecm,height=5truecm]{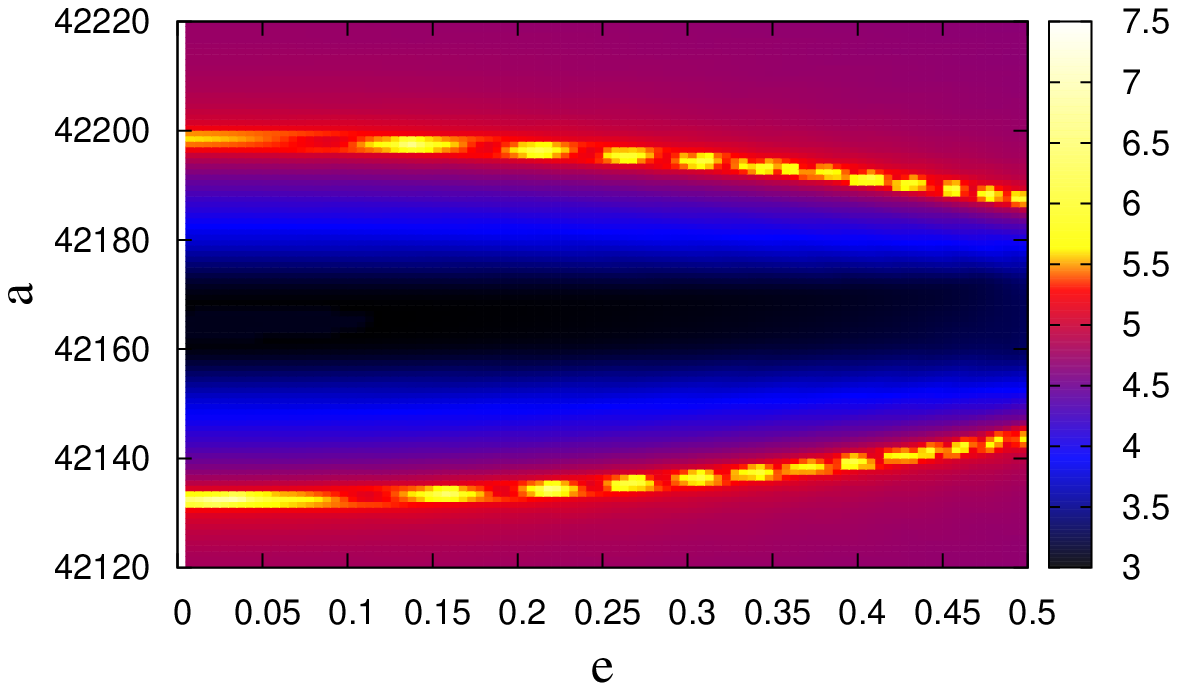}\\
\vglue-0.8cm
\includegraphics[width=6truecm,height=5truecm]{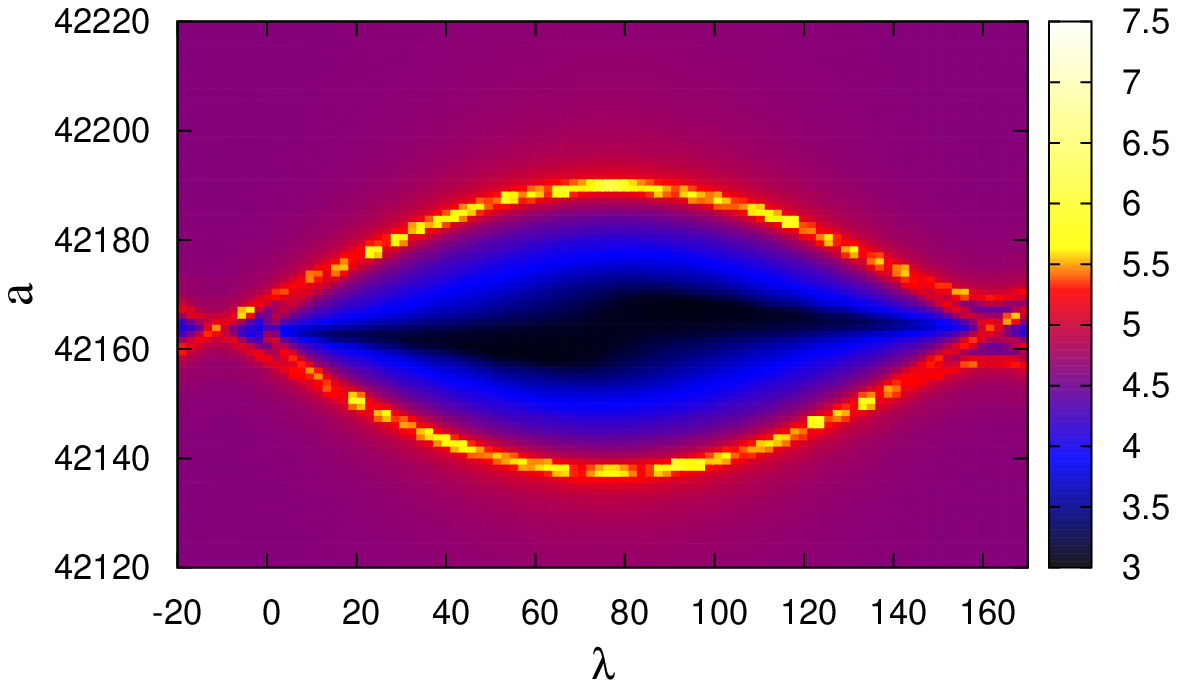}
\includegraphics[width=6truecm,height=5truecm]{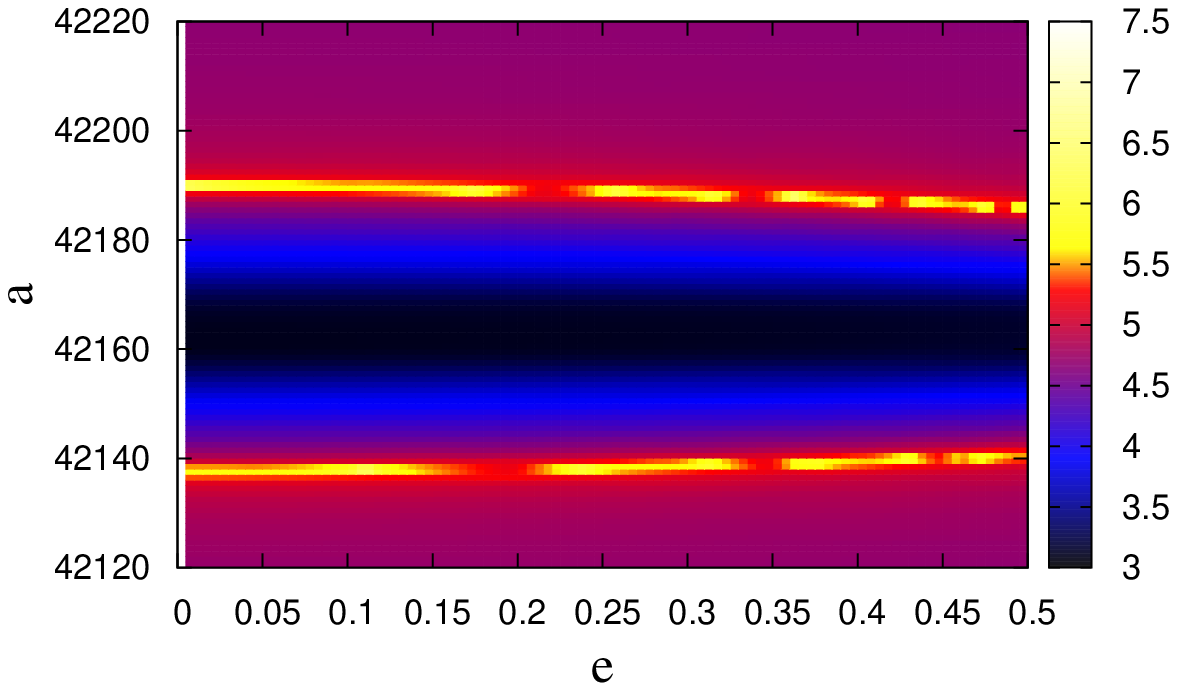}
\vglue0.3cm
\caption{FLI (using Hamilton's equations) for the GEO 1:1 resonance,  under the effects of
all harmonics up to degree and order $n=m=4$, for $e=0.005$, $\omega=0$, $\Omega=0$ and $i=30^o$ in the
upper left panel,  $i=60^o$ in the bottom left panel.
The right panels provide FLI for $\lambda=75.07^o$, $i=30^o$ (top right) and $i=60^o$ (bottom right) in the $(a,e)$ plane.  }
\label{fligeo_i=30}
\end{figure}

\begin{figure}
\centering
\vglue-0.9cm
\includegraphics[width=6truecm,height=5truecm]{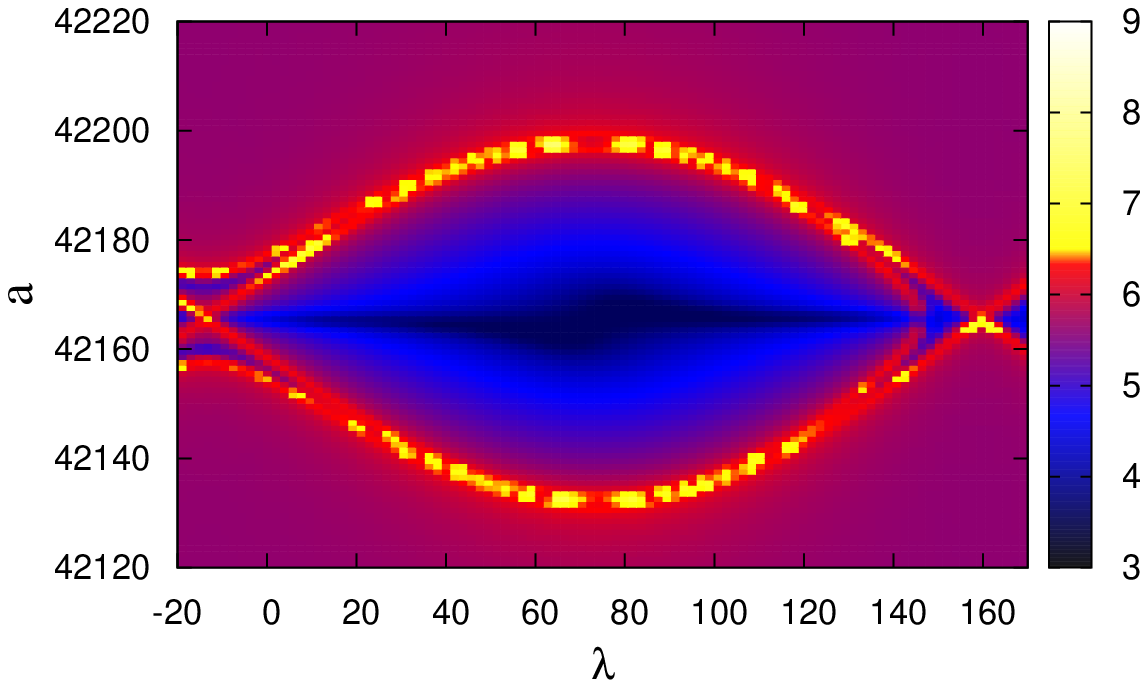}
\includegraphics[width=6truecm,height=5truecm]{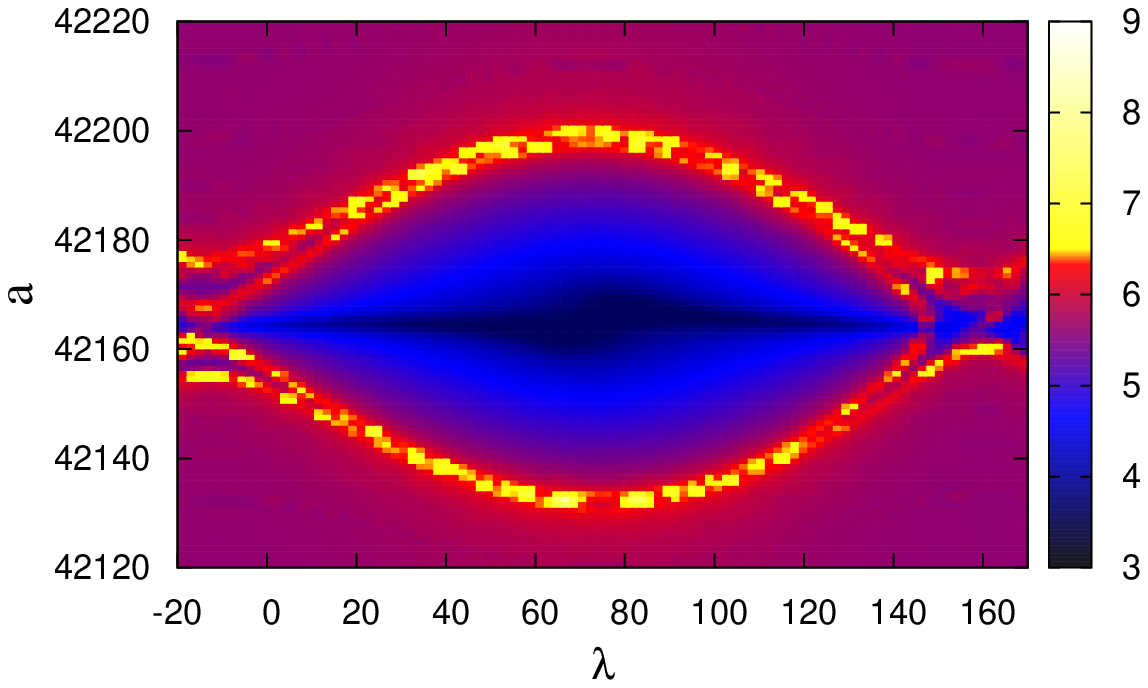}
\vglue0.3cm
\caption{FLI (using Cartesian equations) for the GEO 1:1 resonance for $e=0.005$, $i=30^o$, $\omega=0$, $\Omega=0$,
under all harmonics up to degree and order three (left panel),
all harmonics up to degree and order three + Moon + Sun+ SRP with $A/m=0.1$ (right panel). }
\label{fligeo_cartesian}
\end{figure}

We compute a grid of $100\times 100$ points of
the $\lambda - a$ plane, where the stroboscopic mean node $\lambda$ ranges in the interval $[-220^o,180^o]$, while the semimajor axis $a$ spans an area around the geostationary value corresponding to $a=42\,164.1696$ km.
Figure~\ref{fligeo} shows the FLI values for $e=0.005$, $i=0$, $\omega=0$, $\Omega=0$ and different approximations
from second to fourth order degree harmonics (we used as integration algorithm a 4--th order Runge--Kutta method).
We also add a panel in the $(e,a)$ plane, which provides the amplitude of the libration islands as the eccentricity increases.

For the dynamical model based just on the perturbing harmonics $J_2$ and
$J_{22}$ (Figure~\ref{fligeo}, upper left panel), pendulum like plots are obtained.
The semimajor axis $a$ and the stroboscopic mean node $\lambda$ librate or circulate; the separatrix divides the phase-space in regions corresponding to libration or circulation.
Including higher order harmonics, up to degree and order three and four,
the pattern of the resonance slightly changes, revealing the appearance of more complex orbits
(see Figure~\ref{fligeo}, upper right and bottom left panels).
In particular, near the separatrix some additional curves are visible and, moreover,
the libration zones lose their symmetry.
These complex orbits are not a consequence of the resonance splitting, as in the case of the 2:1 resonance
(see Section~\ref{sec:cart21}), but they are rather due to the interaction between $\mathfrak{T}_1$ and $\mathfrak{T}_3$.
Since the argument of $\mathfrak{T}_1$ is $2(\lambda- \lambda_{22})$, while the argument of $\mathfrak{T}_3$
is $\lambda -\lambda_{31}$, we get an intricate dynamical behavior clearly depicted in Figure ~\ref{fligeo}, upper right panel.
The fourth degree  harmonic terms play a negligible role (compare with Figure~\ref{fligeo}, bottom left).
The pattern shown on the first half of each plot of Figure~\ref{fligeo}, that is for $\lambda \in [-220^o,-20^o]$,
could be viewed as a mirror reflection of the one appearing in the other half. Since all figures obtained for
the $1:1$ resonance have this feature, in the following we plot the FLI values just for  $\lambda $ in the interval $[-20^o,170^o]$.
Increasing the eccentricity we have a decrease of the amplitude of the librational region around the 1:1 resonance
(Figure~\ref{fligeo}, bottom right panel). This effect
is common to all cases shown in the present subsection.

Figure~\ref{fligeo_i=30} shows the FLI values in the spatial case ($i=30^o$ or $i=60^o$) as a function of $\lambda$ and $a$,
or in the $(e,a)$ plane.
All these figures confirm that $\mathfrak{T}_1$ is dominant. As a consequence, the stable point is located at about $\lambda=75^o$ and since $\mathfrak{T}_1$ is proportional to $(1+\cos i)^2 (1-\frac{5}{2}e^2)$, as already remarked the amplitude of the resonance slightly decreases when the inclination and/or the eccentricity increase.

It is worth mentioning that although $\mathfrak{T}_3$ is dominant for inclinations close to $i=60^o$ and large
eccentricities, the dynamics is still leaded by the terms of order $J_{22}$.  For such cases,
the magnitude of $\mathfrak{T}_1 +\mathfrak{T}_2 $ is larger than the magnitude of any combination of terms having
as argument  $\lambda+const$, and as a consequence the equilibrium points are located  at the same points as for
small eccentricities.

The results obtained by using the Hamiltonian formulation are validated by integrating the Cartesian equations of motion as in
Figure~\ref{fligeo_cartesian}, computed for  $i=30^o$, $\omega=0$, $\Omega=0$, $e=0.005$. We have used as starter
a single step method (a Butcher numerical algorithm), while a multistep predictor--corrector
numerical method (Adams-Bashforth 12 steps and Adams-Moulton 11 steps)
performs most of the propagation. The Adams-Bashforth (predictor) method gives an initial solution, while the Adams-Moulton (corrector) method is
successively used to generate better estimates of the solution. Being a convergent process, provided the step size is small enough,
the algorithm was iteratively applied till the solution was obtained within a fixed tolerance.
 The integration step was set to $h=10$ minutes, while the total time span
was $N=15000$ sidereal days (about 41 years). For each orbit we used a fixed initial tangent vector,
leaving the exploration of the effect of a random choice of this vector to a future work.
The left panel of Figure~\ref{fligeo_cartesian} is obtained by considering all harmonics up to degree and order three,
while in the right panel we add the effect of Moon, Sun and SRP for an object with $A/m=0.1$.
The overall structure is very similar to the integration of the Hamilton's equations
in Figure~\ref{fligeo_i=30}. Both approaches predict the same dynamical behavior; precisely:
the equilibrium points have the same location, the long term behavior of the orbital elements are similar,
in particular the semimajor axis is nearly the same, some chaotic orbits are observed near the separatrix.

From the analysis provided in this subsection, we conclude that a dynamical model based just on the perturbing harmonics $J_2$ and $J_{22}$ yields the essential features of the dynamics providing pendulum like plots, while the inclusion of higher order harmonics in the  model provokes a small interaction between all harmonics (tesseral, sectorial and zonal), revealed by the appearance of additional structures near the separatrix and by the lack of symmetry of the librational zones.

For small and moderate inclinations and eccentricities, the effects of Sun, Moon and solar radiation pressure with a small $A/m$ parameter (Figure~\ref{fligeo_cartesian}, right panel) do not change significantly the pattern of the resonance, a fact pointed out also in~\cite{BWM} and \cite{VDLC}.

\subsection{Location of the equilibrium points for the 1:1 resonance}\label{sec:locationGEO}
In the above subsection an analysis of the 1:1 resonance was presented for various eccentricities and inclinations. In all computations, we considered the initial conditions $\omega=0$ and $\Omega=0$. In this section we analyze the influence of non--zero initial angles $\omega$ and $\Omega$ on the location of the equilibrium points.

Since all resonant arguments in $R_{earth}^{res1:1}$ may be written in the form $m \lambda +j \omega -n \lambda_{nm}$, where $m$, $n \in \mathbb{N}$ and  $j \in \mathbb{Z}$, a non--zero initial $\Omega$ does not influence the location of the  equilibrium points. Moreover, since $\mathfrak{T}_1$ is dominant in almost all regions of the phase space, with its magnitude much greater than the size of any other term, and since its resonant argument (see \eqref{t11_terms}) does not depend on $\omega$, the location of the equilibrium points is not affected by the argument of perigee for almost all inclinations and eccentricities.

\begin{figure}
\centering
\vglue-0.9cm
\includegraphics[width=6truecm,height=5truecm]{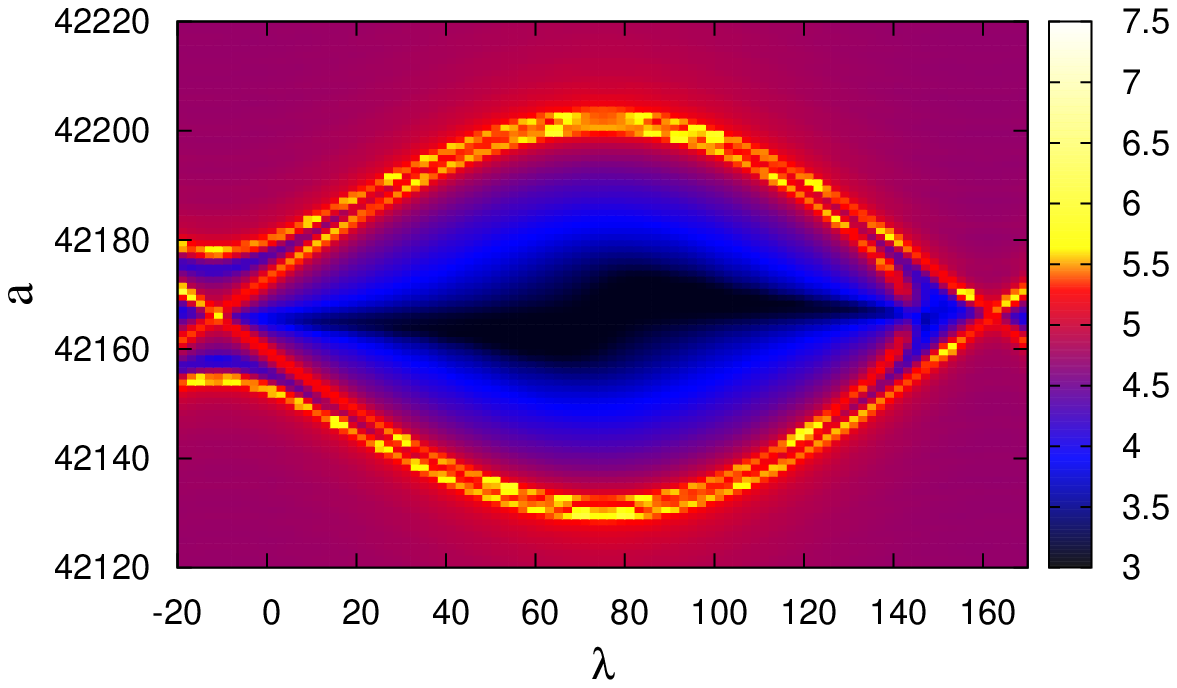}
\includegraphics[width=6truecm,height=5truecm]{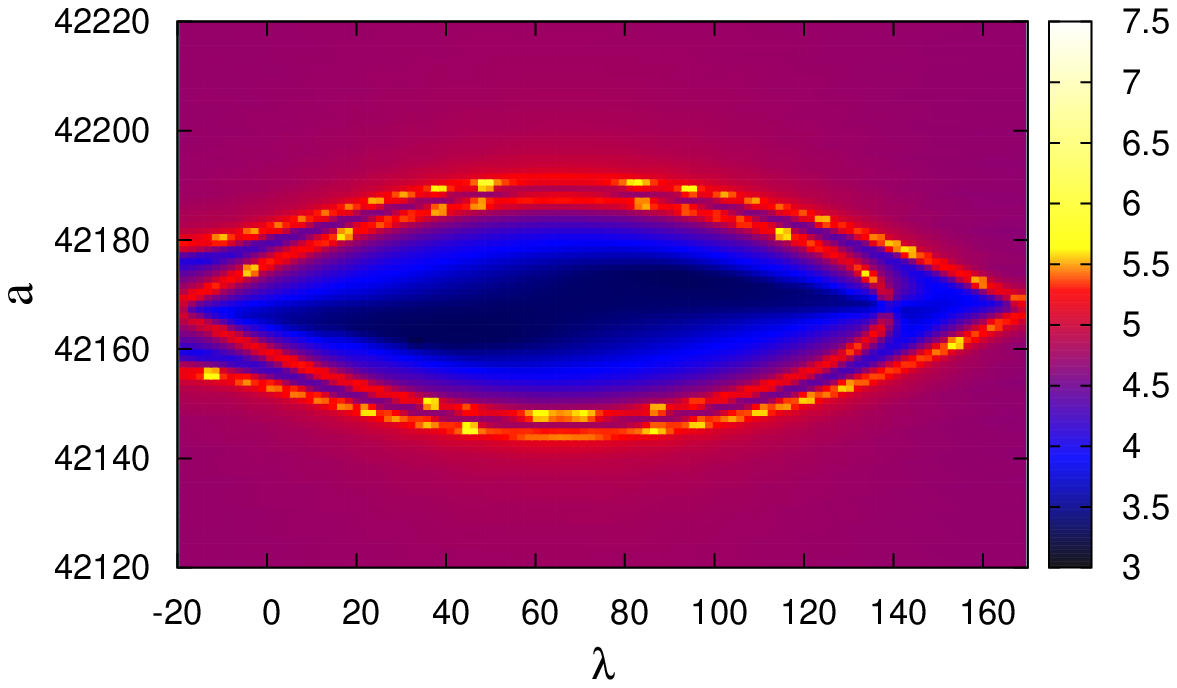}
\vglue0.3cm
\caption{FLI (using Hamilton's equations) for the GEO 1:1 resonance,  under the effects of
all harmonics up to degree and order $n=m=4$, for $i=0$, $\omega=20^o$, $\Omega=0$:  $e=0.005$ (left); $e=0.5$ (right).}
\label{geo_i=0_om=20}
\end{figure}

Figure~\ref{geo_i=0_om=20}, obtained for $\omega=20^o$, $\Omega=0$, $i=0$, $e=0.005$ (left panel) and $e=0.5$ (right panel) provides evidence for
the above claims.
In particular, it is clear that the location of the equilibrium points is not affected by $\omega$.

We immediately recognize that using $\lambda=M-\theta+\omega+\Omega$ in place of $M$, then
$\Omega$ becomes a cyclic variable.
Moreover, in the case of GEO, for a model based just on the perturbations due to
$J_2$ and to $\mathcal{T}_1$, then $\omega$ is also a cyclic variable.

Notice that also in the case of the 2:1 resonance, for the model based on the perturbing terms due to $J_2$ and to the term $t_3$
defined in \equ{t1t2t3},
then $\omega$ is a cyclic variable.

For the models based on $J_2$ and another term, one can easily perform an adapted
change of canonical variables in such a way that $\omega$ becomes cyclic.
This is the reason for having that $\Omega$ has no influence on the location of the equilibrium points and that
the location of the equilibrium points depends on $\omega$, provided that there
are two or three terms in the Hamiltonian which are comparable in magnitude.

\subsection{Cartography of the 2:1 resonance}\label{sec:cart21}
In this section we consider the 2:1 resonance and we perform an analysis similar to that provided for the 1:1 resonance in
Section~\ref{sec:cart11}.

\begin{figure}[h]
\centering
\vglue-0.9cm
\includegraphics[width=6truecm,height=5truecm]{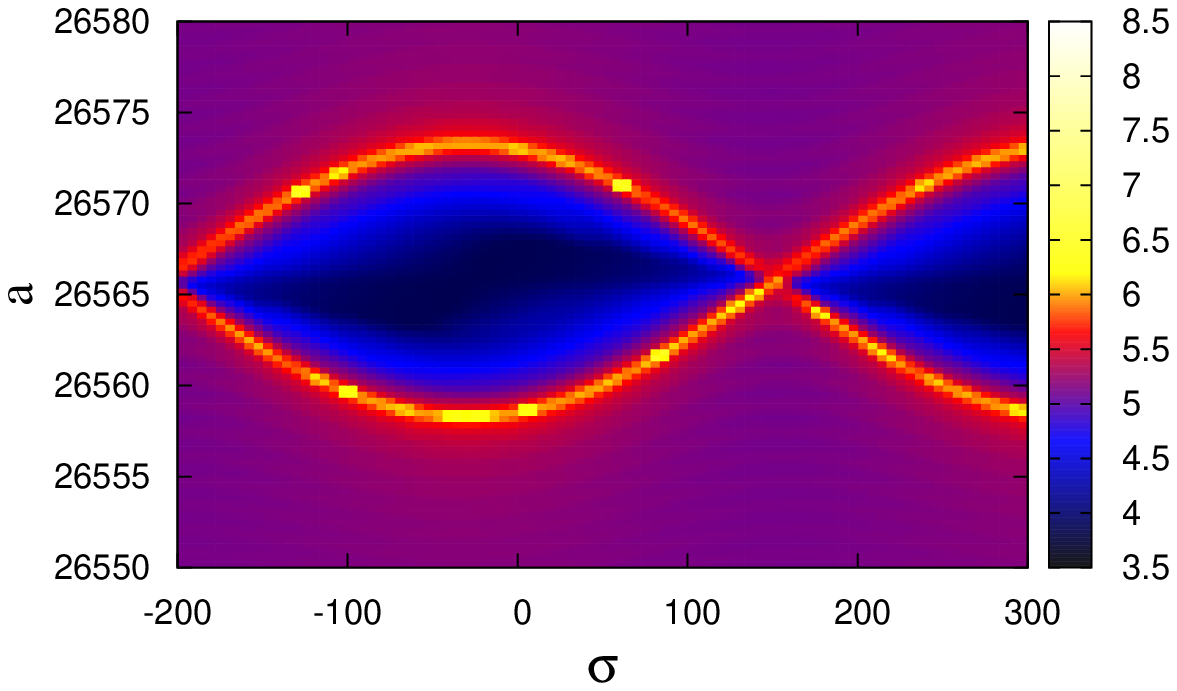}
\includegraphics[width=6truecm,height=5truecm]{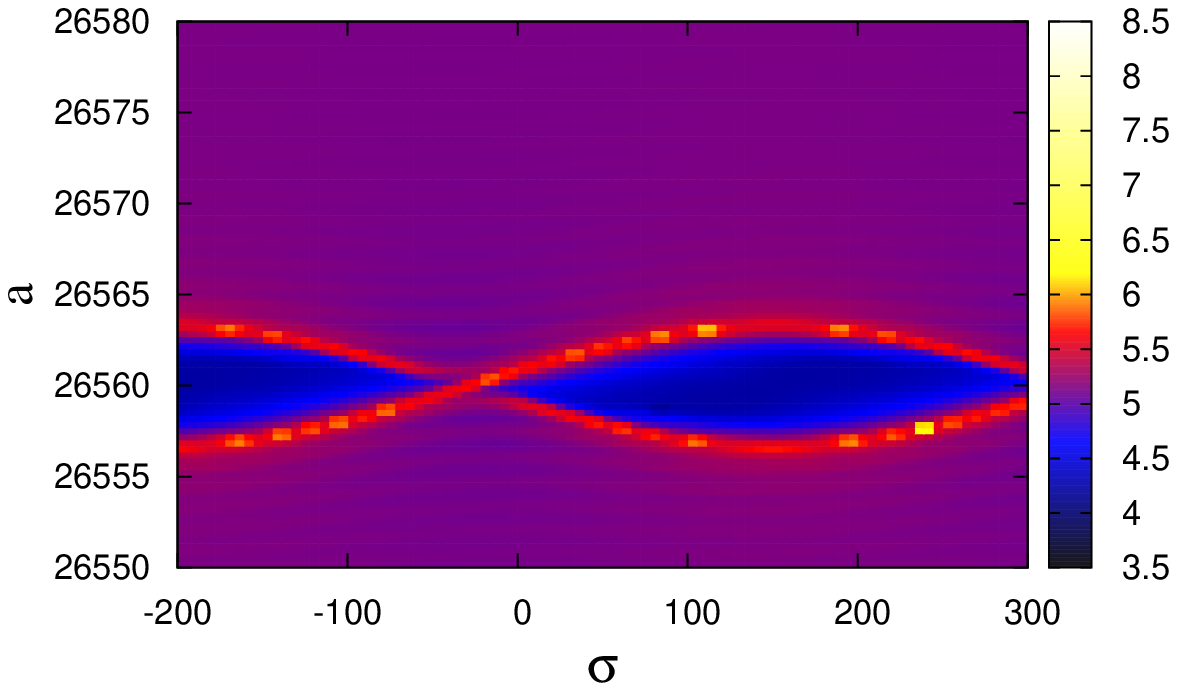}\\
\vglue-0.8cm
\includegraphics[width=6truecm,height=5truecm]{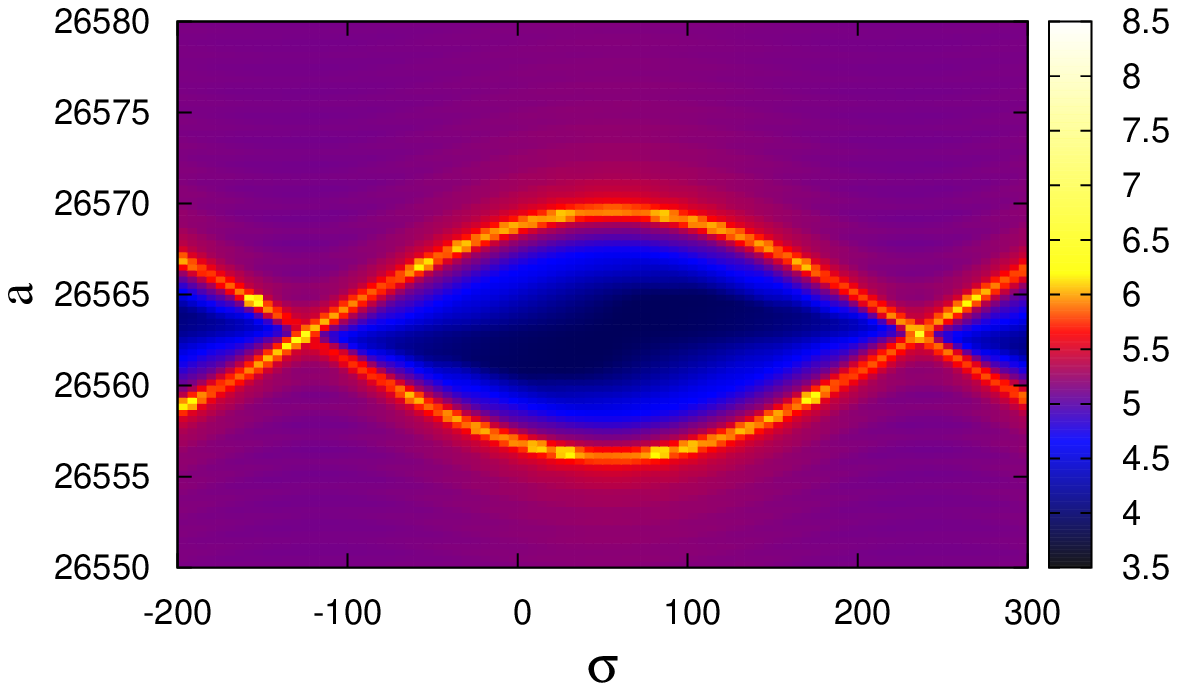}
\includegraphics[width=6truecm,height=5truecm]{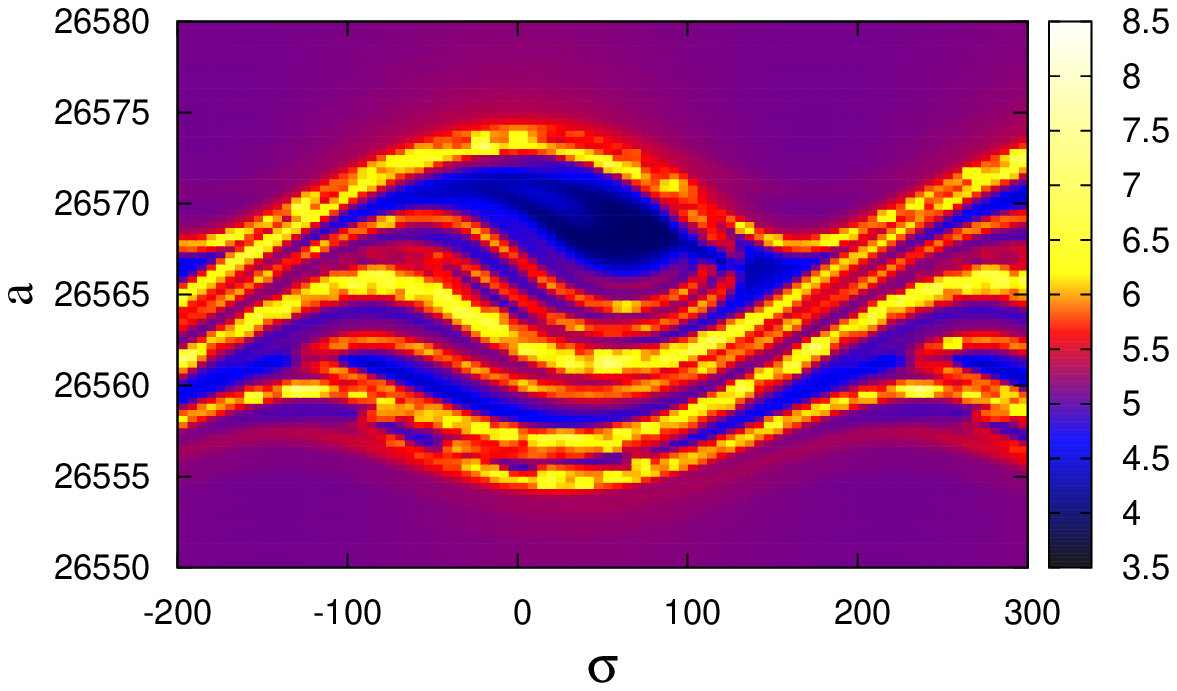}
\vglue0.3cm
\caption{FLI for the toy--model \eqref{toy_ham}, for $e=0.1$, $i=20^o$,
$\omega=0$, $\Omega=0$, under various effects:  $J_2+t_1$ (top left); $J_2+t_2$ (top right);  $J_2+t_3$ (bottom left);
 $J_2+t_1+t_2+t_3$ (bottom right).}
\label{toy}
\end{figure}

\begin{figure}[h]
\centering
\vglue-0.9cm
\includegraphics[width=6truecm,height=5truecm]{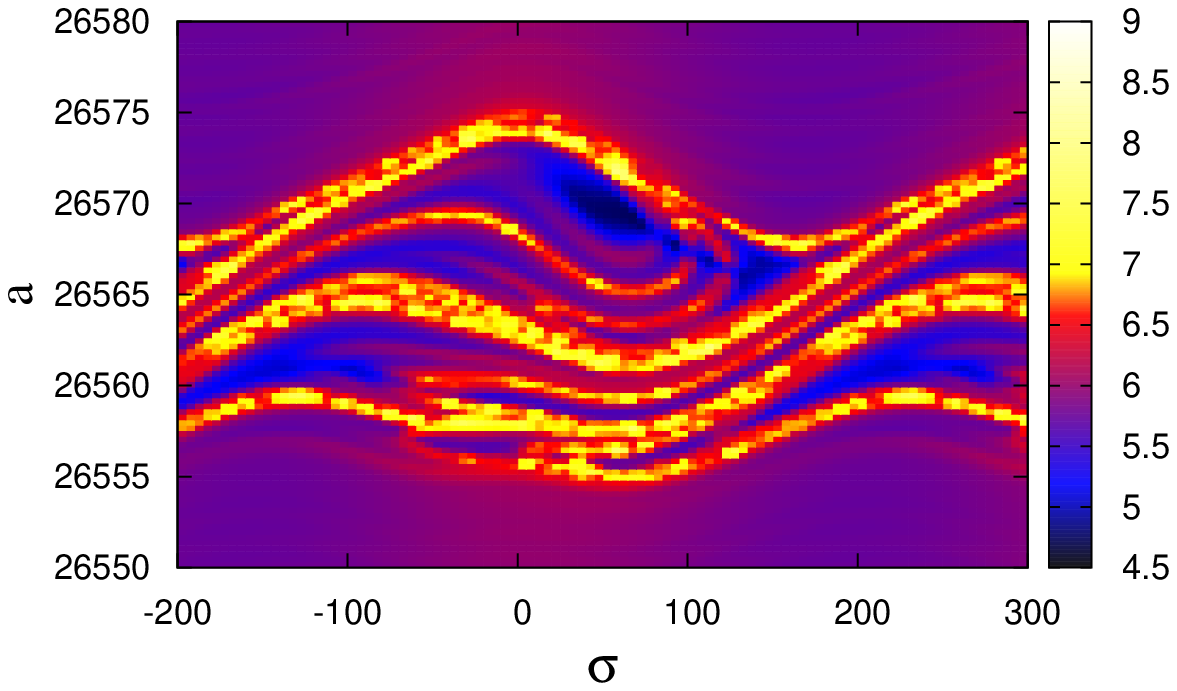}
\includegraphics[width=6truecm,height=5truecm]{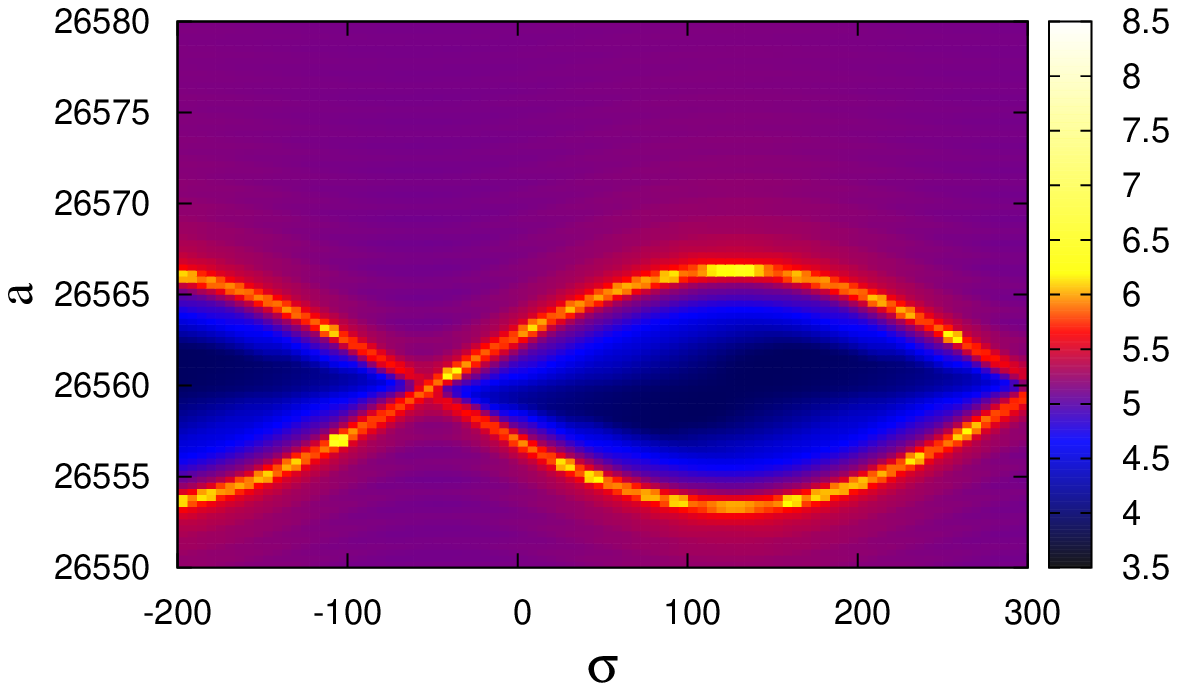}
\vglue0.3cm
\caption{Left: FLI (using Cartesian equations) for the MEO 2:1 resonance for  $i=20^o$, $e=0.1$, $\omega=0$,
$\Omega=0$, under the effects of harmonics up to degree and order 3, Sun,
Moon and solar radiation pressure with $A/m=0.1$.
Right: FLI for the toy--model \eqref{toy_ham}, for $e=0.1$, $i=63.4^o$,
$\omega=0$, $\Omega=0$ under the effects of
 $J_2+t_1+t_2+t_3$.}
\label{meo_i=20_e=0.1}
\end{figure}

For the $2:1$ resonance a phenomenon of superposition of harmonics takes place. To explain this phenomenon and to point out its effects, let us consider the following toy model described by the Hamiltonian
\begin{equation}\label{toy_ham}
\mathcal{H}_{toy}=-\frac{\mu_E^2}{2L^2} +\widetilde R_{earth}^{sec}+t_1+t_2+t_3\ ,
\end{equation}
where we approximate $R_{earth}^{sec}$ by
$$
\widetilde R_{earth}^{sec}\cong\frac{\mu_E R^2_E J_{2}}{a^3}
\Bigl(\frac{3}{4} \sin^2 i -\frac{1}{2}\Bigr) (1-e^2)^{-3/2}\ ,
$$
that is, we consider just the influence of the $J_2$ harmonic, and where $t_1$, $t_2$, $t_3$ are defined in (\ref{t1t2t3}).
For simplicity, we consider these terms up to second order in eccentricity, namely we take
\beqano
t_1&=& \frac{\mu_E R_E^2 J_{22}}{a^3} \Bigl\{\frac{3}{4} (1+\cos i)^2 \Bigl( -\frac{e}{2}\Bigr) \cos (\sigma +\omega -2\lambda_{22})\Bigr\}\nonumber\\
t_2&=&\frac{\mu_E R_E^2 J_{22}}{a^3} \Bigl\{\frac{3}{2} \sin^2 i
\Bigl(\frac{3}{2}e \Bigr)  \cos (\sigma -\omega-2\lambda_{22}) \Bigr\} \nonumber\\
t_3&=& \frac{\mu_E R_E^3 J_{32}}{a^4} \Bigl\{ \frac{15}{8} \sin i
(1-2\cos i-3 \cos^2 i) \Bigl(1+2e^2 \Bigr) \sin(\sigma -2\lambda_{32})
\Bigr\}\ ,
\eeqano
where $\sigma=2 \lambda$ with $\lambda$ as in \equ{lambda_definition_2:1}.

\begin{figure}[h]
\centering
\vglue-0.9cm
\includegraphics[width=6truecm,height=5truecm]{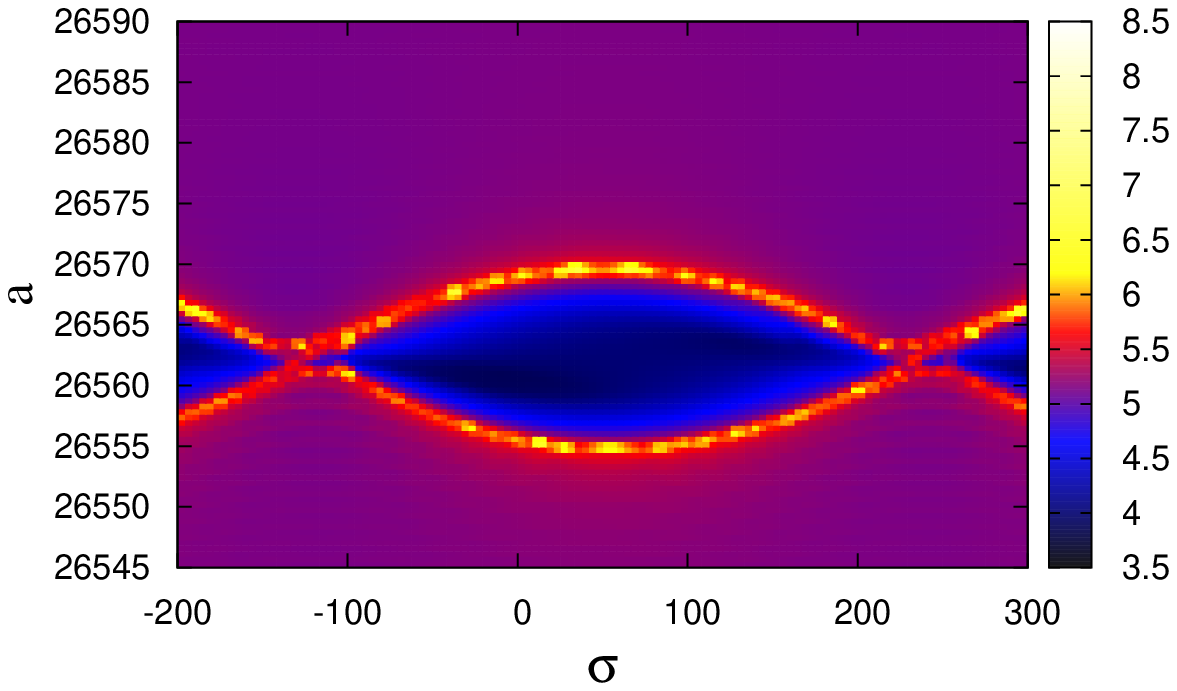}
\includegraphics[width=6truecm,height=5truecm]{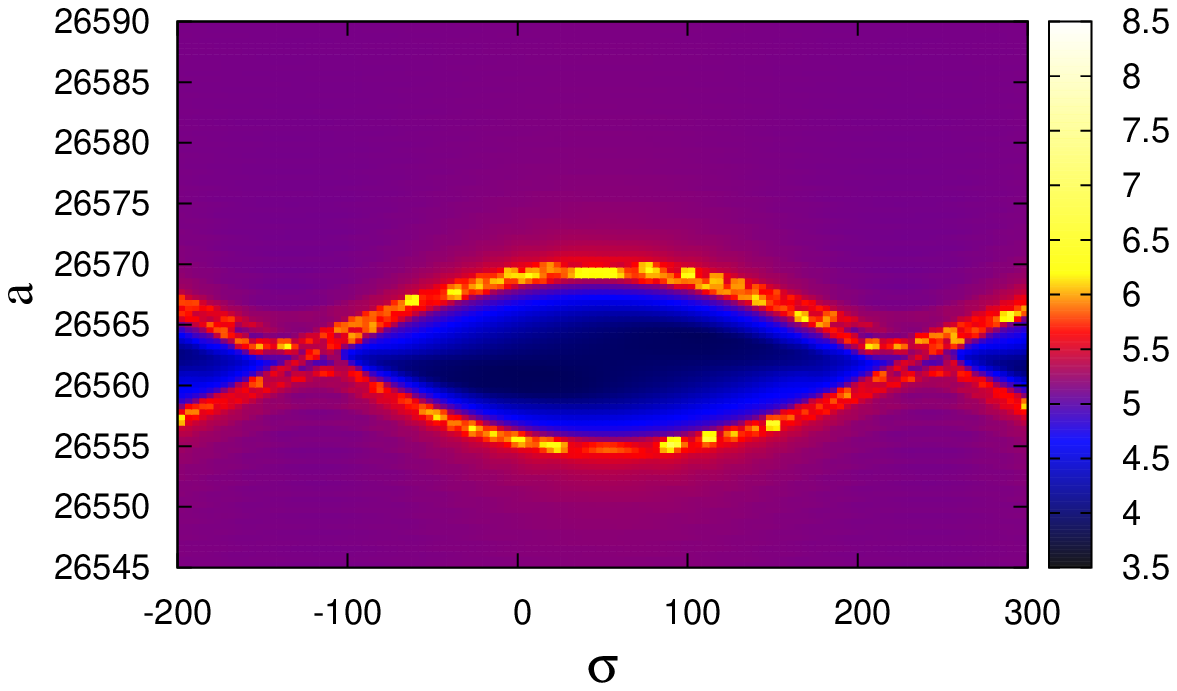}\\
\vglue-0.8cm
\includegraphics[width=6truecm,height=5truecm]{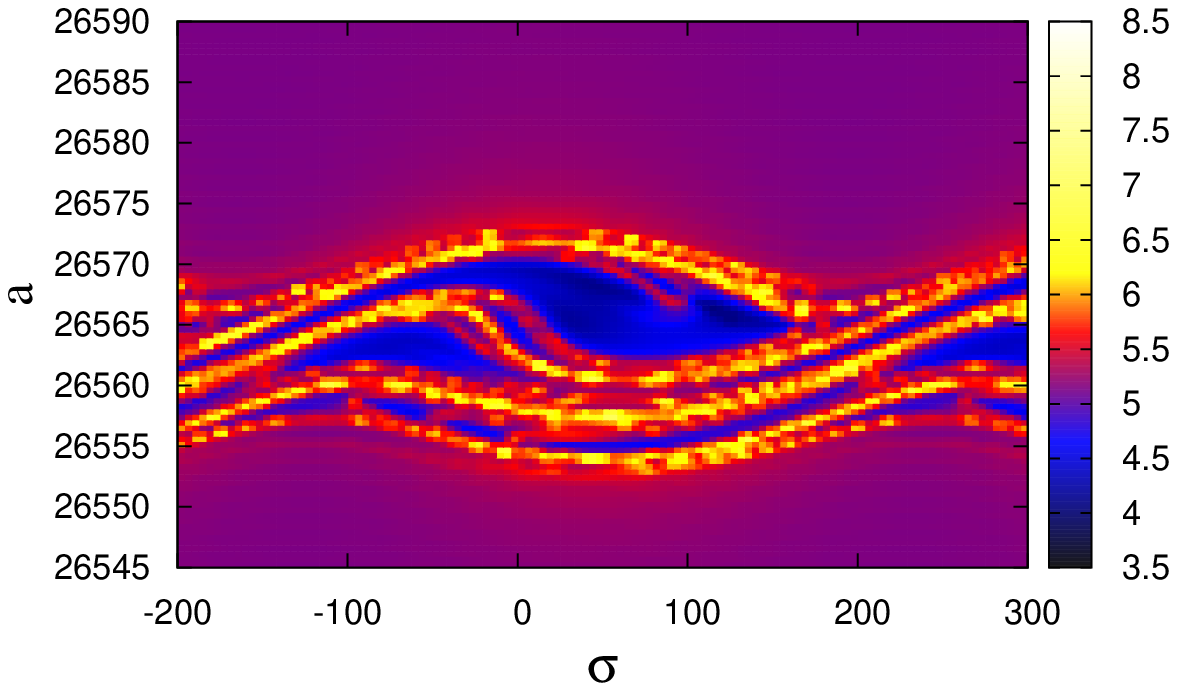}
\includegraphics[width=6truecm,height=5truecm]{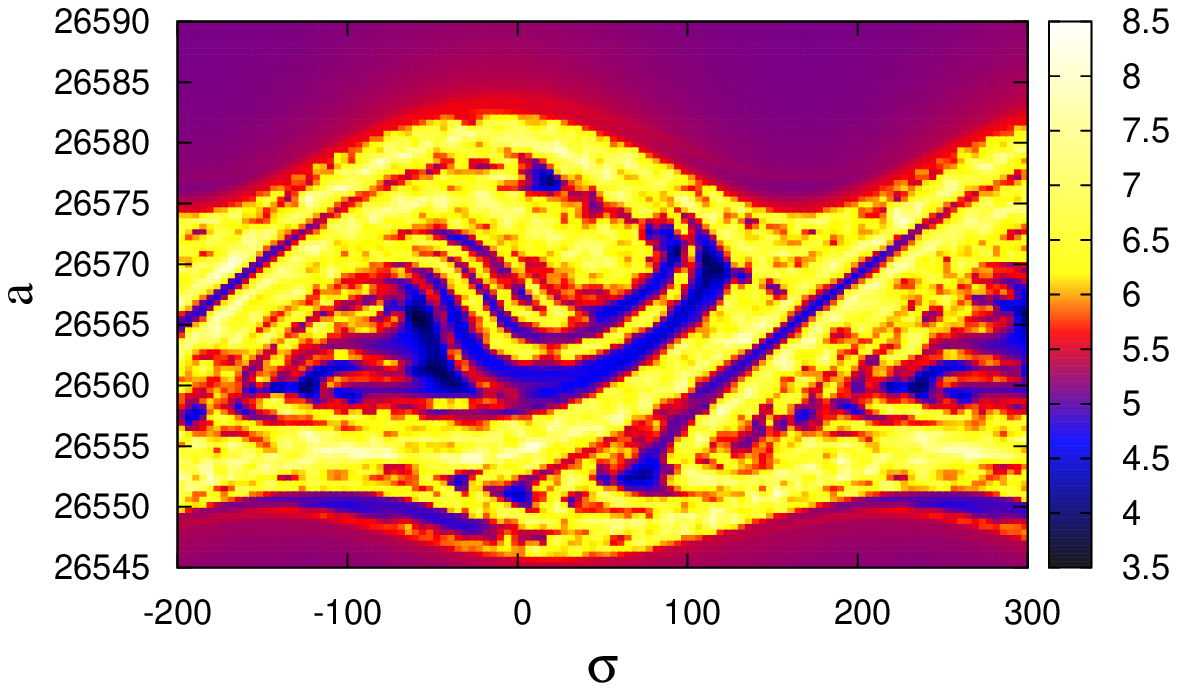}
\vglue0.3cm
\caption{FLI (using Hamilton's equations) for the MEO 2:1 resonance,  under the effects of
all harmonics up to degree and order $n=m=4$, for $i=30^o$, $\omega=0$, $\Omega=0$:  $e=0.005$ (top left);  $e=0.01$ (top right);  $e=0.1$ (bottom left); $e=0.5$ (bottom right).}
\label{res21_i=30}
\end{figure}

It is worth mentioning that, for small and moderate eccentricities, this simple dynamical model yields the essential features of the dynamics inside the 2:1 resonance. The other terms of $R_{earth}^{sec}$ and $R_{earth}^{res2:1}$ have a secondary role.

The canonical variables in the non--autonomous Hamiltonian \eqref{toy_ham} are the action--angle Delaunay variables
$(L,G,H,M,\omega,\Omega)$, which are related to the orbital elements by the relations \eqref{Delaunay_var}. An autonomous Hamiltonian $\widetilde{\mathcal{H}}_{toy}$ can be obtained from $\mathcal{H}_{toy}$ by introducing an artificial momentum $\Theta $
conjugated to $\theta$. Precisely, assuming that $\dot\theta=n_\theta$, we consider the Hamiltonian
$$\widetilde{\mathcal{H}}_{toy}=\mathcal{H}_{toy}+n_\theta \Theta\ ,$$
where $\dot\theta={{\partial\widetilde{\mathcal{H}}_{toy}}\over {\partial\Theta}}$, $\dot\Theta=-
{{\partial\widetilde{\mathcal{H}}_{toy}}\over {\partial\theta}}$.
Let us consider the symplectic canonical transformation:
$$(L,G,H, \Theta, M,\omega,\Omega, \theta) \longrightarrow (L^\prime,G^\prime,H^\prime, \Theta^\prime, \sigma,\omega,\Omega, \theta)\ ,$$
which is defined by
\begin{equation} \label{new_var}
L^\prime=L\ , \quad G^\prime=G-L\ , \quad H^\prime=H-2L\ , \quad \Theta^\prime=\Theta +2L\ ,\quad \sigma=M-2\theta+\omega+2\Omega\ ,
\end{equation}
while $\omega$, $\Omega$, $\theta$ are kept unaltered.
In terms of the new variables, the autonomous Hamiltonian $\widetilde{\mathcal{K}}_{toy}$ is given by
$$
\widetilde{\mathcal{K}}_{toy}=-\frac{\mu_E^2}{2L'^2}-2 n_\theta L'+n_\theta\Theta' +\widetilde R_{earth}^{sec}+t_1+t_2+t_3\ ,
$$
from which we get the non--autonomous Hamiltonian
\begin{equation}\label{K_ham}
\mathcal{K}_{toy}=-\frac{\mu_E^2}{2L'^2}-2 n_\theta L'+\widetilde R_{earth}^{sec}+t_1+t_2+t_3\ .
\end{equation}
The new action--angle variables $(L,G^\prime,H^\prime,\sigma,\omega,\Omega)$ are related to the orbital elements through the relations \eqref{Delaunay_var}, \eqref{lambda_definition_2:1}, \eqref{new_var}.

Let us underline one of the main effects induced by the $J_2$ harmonic (see also \cite{chao}).
Setting $t_1=t_2=t_3=0$ in \eqref{K_ham}, it is easy to show (see
for example~\cite{Kaula})  that $J_2$ provokes a secular regression of the
orbital node and a precession of the perigee, that is
\beq{J2_effect}
\dot\Omega \simeq -{3\over 2}n^\ast J_2({R_E\over{a(1-e^2)}})^2\cos i\ ,\qquad
\dot\omega \simeq {3\over 2}n^\ast J_2({R_E\over{a(1-e^2)}})^2(2-{5\over 2}\sin^2 i)\ ,
\eeq
where $n^\ast$ denotes the mean motion.
Let us remark that for $i = 63.4^0$, called the \sl critical inclination, \rm we have $2-{5\over 2}\sin^2 i=0$ and therefore $\omega$ is constant, giving place to the so--called \sl frozen orbits \rm (\cite{chao}).

If $i \neq 63.4^0$, then $\dot{\omega} \neq 0$ and each of the angles $\sigma$, $\sigma+\omega$, $\sigma-\omega$
appearing in each term of \equ{t1t2t3} will have
zero derivatives at different locations. Therefore, the gravitational resonance splits
into a triplet of resonances, the separation between them being of the order of a
few kilometers (precisely a value between 0-9 km, the exact value depending on the
inclination and the eccentricity).

In Figure~\ref{toy}, obtained for the initial conditions $e=0.1$, $i=20^o$, $\omega=0$, $\Omega=0$, we represent the FLI values for the toy-models taking into account as perturbations: $J_2+t_1$ (top left), $J_2+t_2$ (top right), $J_2+t_3$ (bottom left), respectively. Each toy--model yields a pendulum like plot.
In particular, for $J_2+t_1$ the stable point is located at  $\sigma=2 \lambda_{22}-180^o\simeq -30^o$ and $a=26565.8$ km.
For $J_2+t_2$ the stable point is located at $\sigma=2 \lambda_{22}\simeq 150^o$ and $a=26559.9$ km, therefore there is a shift of about $6$ km in semimajor axis, if we compare it with the plot obtained for the previous toy--model (Figure~\ref{toy}, top left).
In the case of $J_2+t_3$ the stable point is located at  $\sigma=2\lambda_{32}-90^o \simeq 55^o$ and $a=26\,562.8$ km.

A different behavior is found for
the model taking into account all perturbing terms, that is $J_2+t_1+t_2+t_3$ (Figure~\ref{toy}, bottom right); all harmonics interact, leading to a complex dynamics with a complicated interplay of regular and chaotic motions.

Figure~\ref{meo_i=20_e=0.1} left is obtained by integrating the Cartesian equations and taking into account all harmonics up to degree and order 3, Sun, Moon and solar radiation pressure with $A/m=0.1$; the plot confirms what was anticipated and described by the above toy--model, namely a complex dynamics due to a strong effect provoked by the interaction of all harmonics.

Roughly speaking, as long as the magnitude of one of the terms $t_1$, $t_2$ or $t_3$ is much greater than the magnitude of the other two, then a pendulum--like pattern is obtained.
On the contrary, when two (or all three) terms $t_1$, $t_2$, $t_3$ are comparable in magnitude and $i\neq 63.4^o$, then all harmonics superpose, leading to chaotic motions as in Figure~\ref{meo_i=20_e=0.1} (left panel).
There exists an exception, which corresponds to the critical inclination $i=63.4^o$ (see Figure~\ref{meo_i=20_e=0.1}, right). In this case, $\omega$ is constant and $\sigma$, $\sigma+\omega$, $\sigma-\omega$ will have zero derivatives at the same location, precisely at about $a=26\,560$ km.

Figure~\ref{res21_i=30} is obtained including all harmonics up to degree and order $4$; the different panels illustrate graphically the phenomenon of superposition of harmonics, for various eccentricities and
for $i=30^o$ (similar results are obtained for other inclinations).

\begin{figure}[h]
\centering
\vglue-0.9cm
\includegraphics[width=6truecm,height=5truecm]{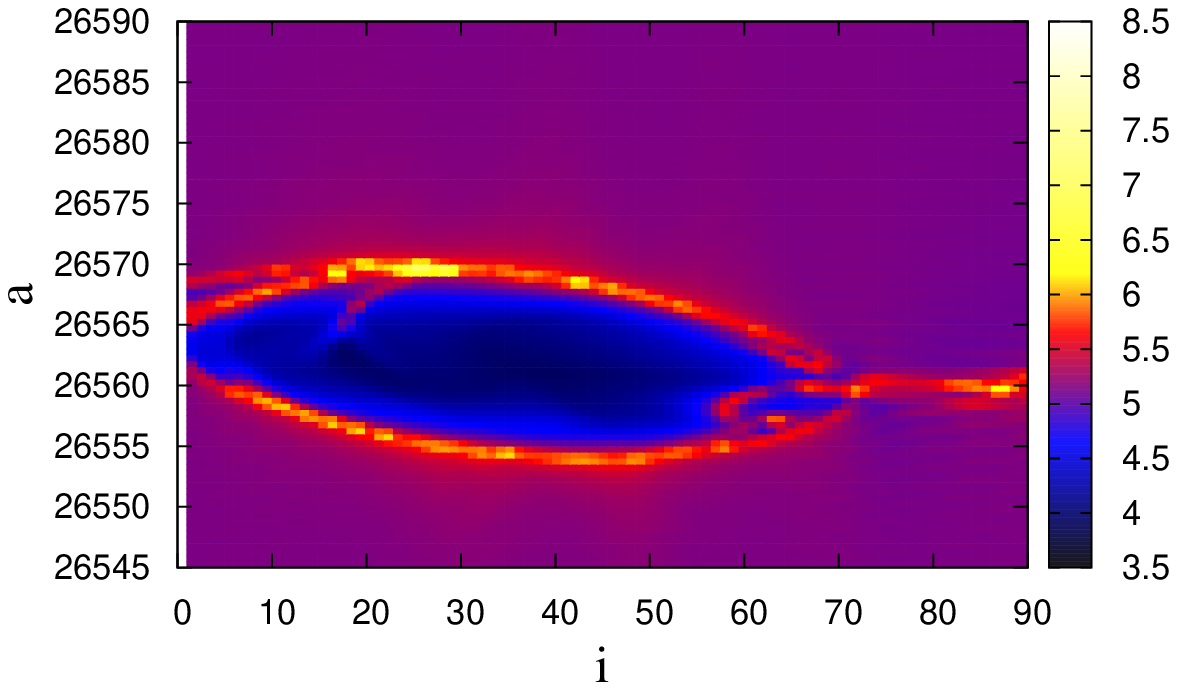}
\includegraphics[width=6truecm,height=5truecm]{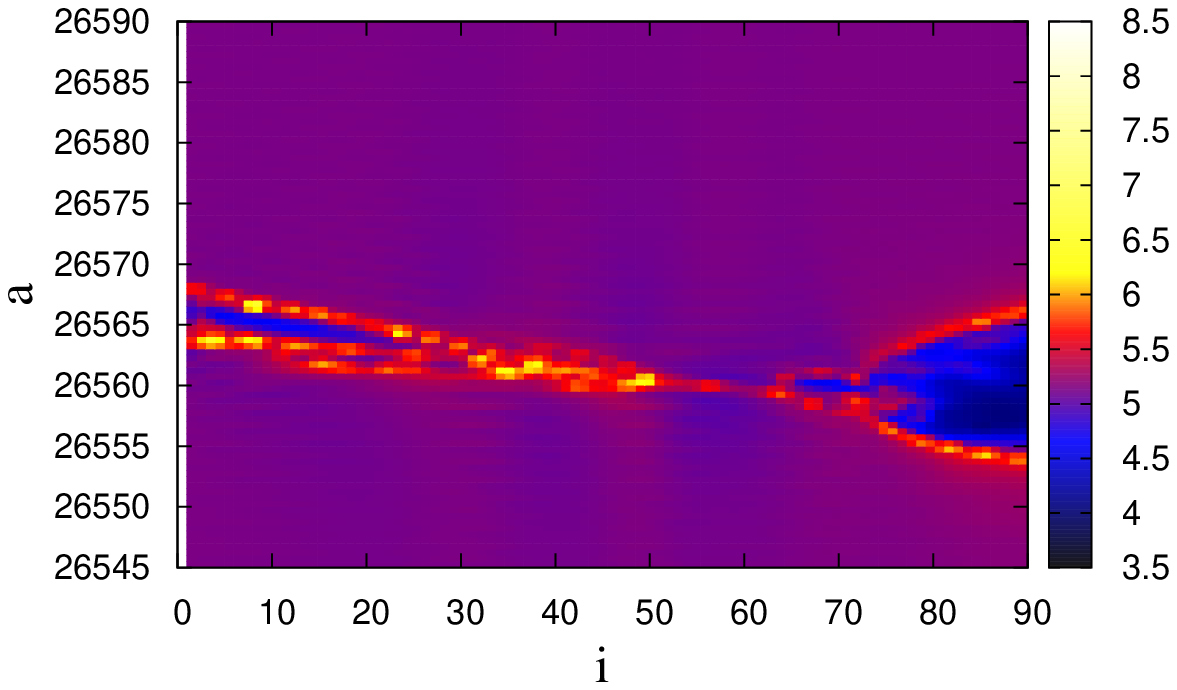}\\
\vglue-0.8cm
\includegraphics[width=6truecm,height=5truecm]{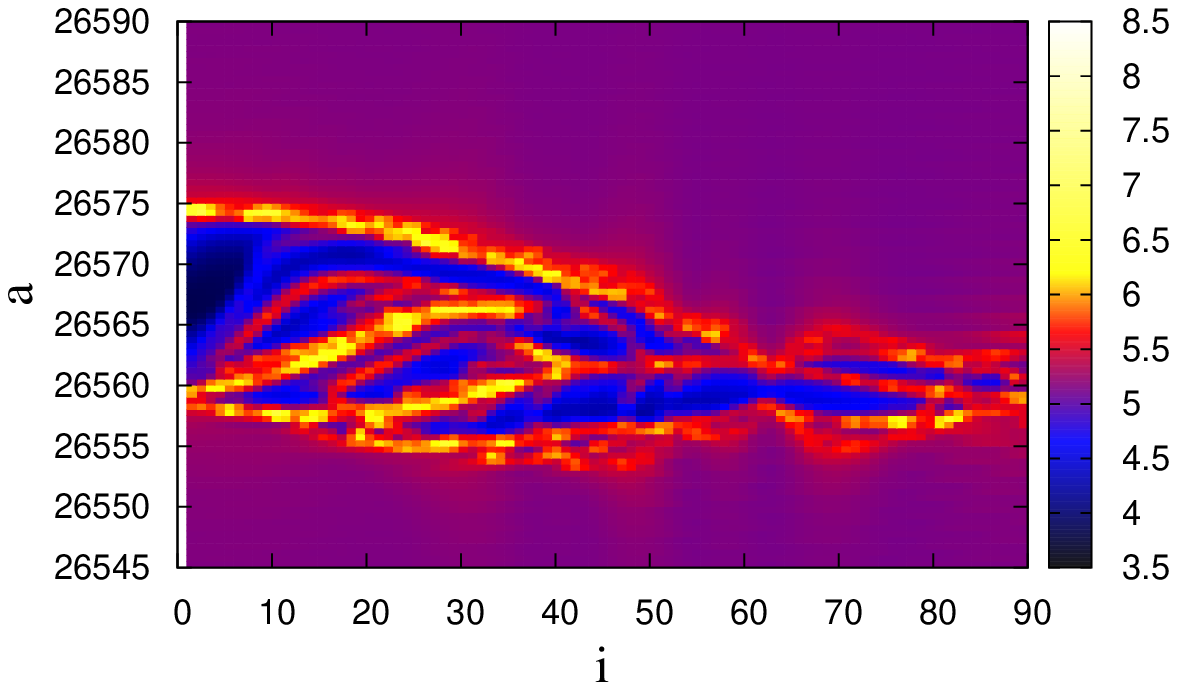}
\includegraphics[width=6truecm,height=5truecm]{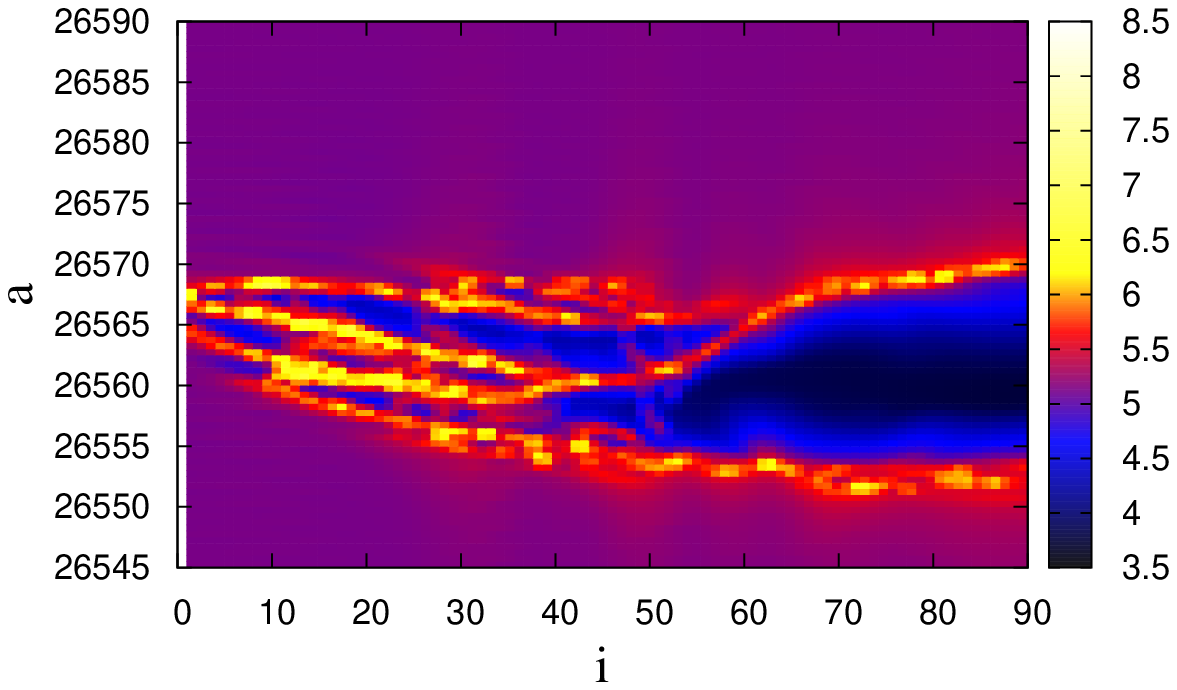}\\
\vglue-0.8cm
\includegraphics[width=6truecm,height=5truecm]{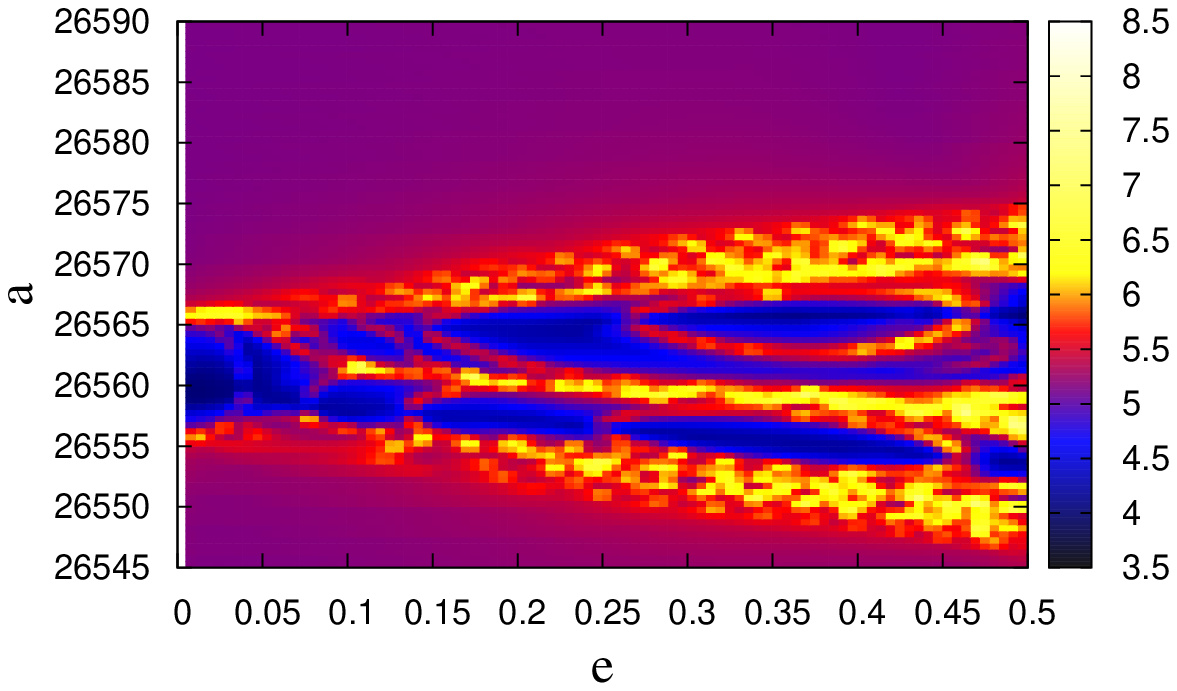}
\includegraphics[width=6truecm,height=5truecm]{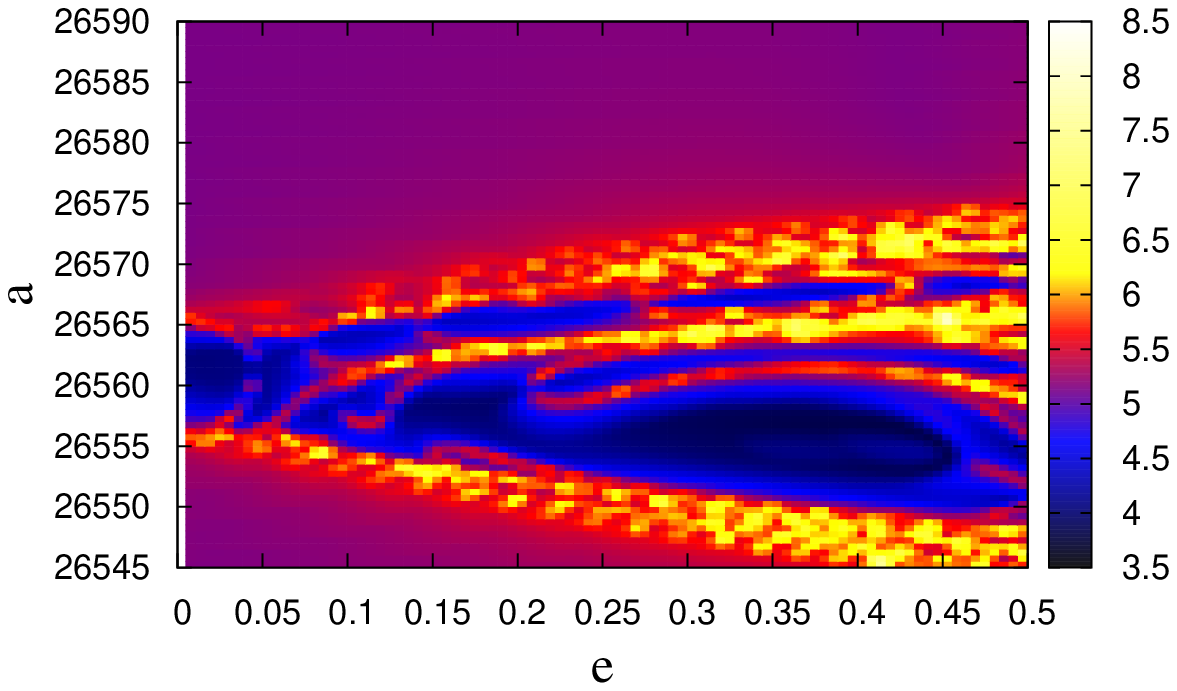}
\vglue0.3cm
\caption{FLI (using Hamilton's equations) for the MEO 2:1 resonance under the effects of all harmonics up to degree and order $n=m=4$.
Upper left and right panels: FLI in $(i,a)$ for $\omega=0$, $\Omega=0$:  $e=0.005$, $\sigma=55^o$ (left); $e=0.005$, $\sigma=235^o$ (right).
Middle and bottom panels provide the FLI for $\omega=0$, $\Omega=0$, $\sigma=-30^o$ in the left panels and $\sigma=150^o$
in the right panels. The middle panels provide the FLI for $e=0.1$ in $(i,a)$, while the bottom panels are for $i=50^o$  in
$(e,a)$.}
\label{res21_e=0005_e=001}
\end{figure}

Since $t_1$ and $t_2$ are proportional to the eccentricity, $t_3$ is proportional to $\sin i(1-2 \cos i -3 \cos^2 i)$ and  $t_2$ is proportional to  $\sin^2 i$, each term will dominate in a specific region of the phase space (see Figure~\ref{big11} right). As a consequence, the location of resonant island's centers as well as the widths of the resonances depend on the values of the eccentricity and of the inclination. For example, in Figure~\ref{res21_i=30} (top left and top right), since $t_3$ is dominant for small eccentricities, the center of the resonant island is at about $\sigma=55^o$ and the amplitude of the resonance increases (and then decreases) with the growth of the inclination. An interesting phenomenon occurs for $i=70.53^o$, since the function
 \begin{equation}\label{bif_fun}
f:[0^o,90^o] \longrightarrow \mathbb{R}, \qquad  f(i)=-\sin i ( 1-2 \cos i -3 \cos^2 i),
 \end{equation}
changes its sign at $70.53^o$. Thus, for $i<70.53^o$ the stable point is located at about $\sigma=55^o$ and the unstable one at  about $\sigma=235^o$, while for $i>70.53^o$ the situation is opposite: the hyperbolic point is located at
$\sigma=55^o$ and the stable one at $\sigma=235^o$ (compare with Figure~\ref{res21_e=0005_e=001}). Therefore, at $i=70.53^o$ we have a transcritical  bifurcation point, since a small change in inclination causes the stability of the equilibrium points to change. This critical point is clearly marked as a cusp in the shape displayed in Figure~\ref{big11} right (compare also with Figure~\ref{fig:amplitude} right).

For higher but moderate eccentricities, let us say $e=0.1$, $t_3$ is still dominant for $i\in [25^o, 45^o]$, but the magnitudes of all three terms $t_1$, $t_2$, $t_3$ are comparable. Due to resonance splitting, each component of the triplet gives rise to a resonance at a specific exact location. Usually, these resonances have amplitude greater than the distance which separates them, so that a superposition phenomenon takes place, leading to complex dynamics with both regular and chaotic behaviors (compare with Figure~\ref{res21_i=30}, bottom left).

If the value of the eccentricity is larger, then $t_1$ is dominant for $i<i_0$ and $t_2$ is dominant for $i>i_0$, where $i_0 \in (0^o,\ 90^o)$ is the solution of the equation
\beq{i0}
\frac{3}{4}(1+\cos i)^2 (\frac{e}{2}-\frac{e^3}{16})-\frac{3}{2} \sin^2 i (\frac{3}{2} e+\frac{27}{16}e^3)=0\ .
\eeq
For $e=0.5$, the value of $i_0$ is about $39^o$. This is the reason for having the centers of the islands located somewhere between $\sigma=-50^o$ and $\sigma=0^o$  in the plots obtained for $i<39^o$ (see bottom right in  Figure~\ref{res21_i=30}). Clearly, for $e=0.5$ and non-zero inclinations, many other terms in $R_{earth}^{res2:1}$ grow in magnitude, their interaction leading to chaotic motions.

Figure~\ref{res21_e=0005_e=001} (top panels)
plots the FLI values as a function of inclination and semimajor axis in order to evaluate the width of the resonance for each value of the inclination and to give a hint on the dynamics inside the resonance. For small eccentricities, let us say $e=0.005$,
$t_3$ is dominant and its magnitude is much greater than the magnitude of the other terms. Therefore, we vary $\sigma$ as $\sigma=55^o$ in Figure~\ref{res21_e=0005_e=001} top left, and $\sigma=235^o$ in the top right panel. These plots give an estimate of the amplitude of the resonance as a function of the inclination. For a specific inclination, the amplitude can be determined by measuring the distance between the two points on the separatrix obtained as the intersection of the vertical line corresponding to that specific inclination and the structure visible on the plot. Due to the above described bifurcation phenomenon, Figure~\ref{res21_e=0005_e=001} top left panel should be used for $i \in [0^o, 70.53^o)$, while Figure~\ref{res21_e=0005_e=001} top right panel for $i>70.53^o$. In fact, the way in which the amplitude of the resonance varies with the inclination is described by the function $f$, defined by \eqref{bif_fun}. The function $f$ increases on the interval $i \in [0^o, 34.42^o]$ and then it decreases. At $i=70.53^o$, $f$ changes its sign, so that on the interval $[70.53^o, 90^o]$ the function $-f$ increases. It is striking to compare the behavior of the function $f$ with the pattern followed in Figure~\ref{res21_e=0005_e=001} top panels by the
amplitude of the resonance. Precisely, the amplitude increases in the interval $[0^o, 34.42^o]$, it decreases in the interval $[34.42^o, 70.53^o]$, and then it increases again in the interval $[70.53^o, 90^o]$.
Analyzing carefully Figure~\ref{res21_e=0005_e=001} top panels, we notice that the structure of the left plot is not horizontal,  but rather slightly inclined. That is, increasing the inclination, the location of the equilibrium points is shifted in semimajor axis.
The $J_2$ harmonic is responsible for this behavior, which provokes the secular regression of the orbital node and the precession of perigee. The exact location of the resonance is given by $\dot{\sigma}\equiv\dot{M}-2\dot{\theta}+\dot{\omega}+2 \dot{\Omega} =0$. In view of \eqref{J2_effect}, one can write
$\frac{d(\dot{\omega}+2 \dot{\Omega})}{d i}=\frac{3}{2} n^\ast J_2 \Bigl(\frac{R_E}{a(1-e^2)}\Bigr)^2 (-5 \cos i +2) \sin i < 0$, for  $i\in(0^o, 66.4^o)$. Therefore $\dot{\omega}+2 \dot{\Omega}$ decreases with the inclination within $(0^o, 66.4^o)$ and, as a consequence, the equilibrium points are shifted in semimajor axis, as far as the inclination increases.

Figure~\ref{res21_e=0005_e=001} (middle panels)
is obtained for $e=0.1$ and shows that the roles played by $t_1$ and $t_2$ enhance the complexity of the problem.
Finally, Figure~\ref{res21_e=0005_e=001} (bottom panels) provides the FLI values in the plane $(e,a)$. These plots give an estimate of the width of the resonance for each value of the eccentricity;
we may conclude that for small $e$ the motion is regular, while for moderate and large $e$ the plots show very complex behaviors.

\subsection{Location of the equilibrium points for the 2:1 resonance}
As in in Section~\ref{sec:locationGEO}
we infer that $\Omega$ does not influence the location of the equilibria. In contrast to the $1:1$ resonance,  $\omega$ plays an important role for moderate and high eccentricities. For the $2:1$ resonance there are three leading terms, $t_1$, $t_2$, $t_3$, with  comparable magnitude in most of the phase space. Since their resonant arguments are $\sigma\pm\omega-2\lambda_{22}$ and $\sigma-2\lambda_{32}$, the location of the equilibria and the pattern of the resonances are strongly affected by $\omega$.

\begin{figure}[h]
\centering
\vglue-1.cm
\includegraphics[width=6truecm,height=5truecm]{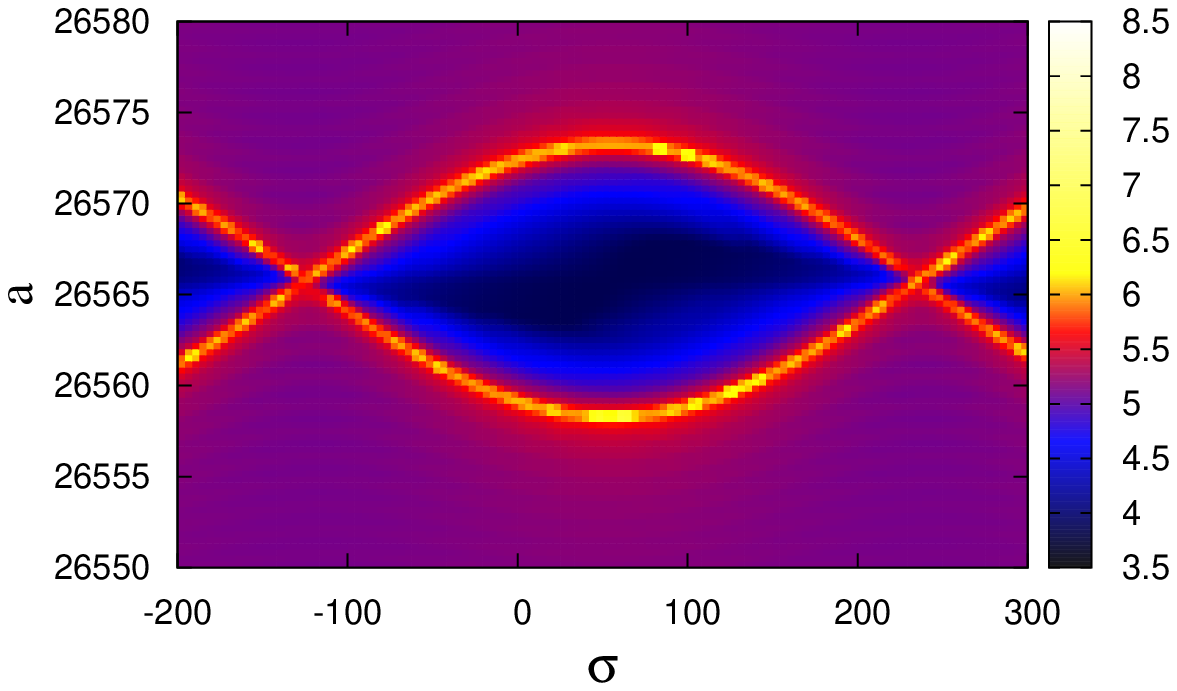}
\includegraphics[width=6truecm,height=5truecm]{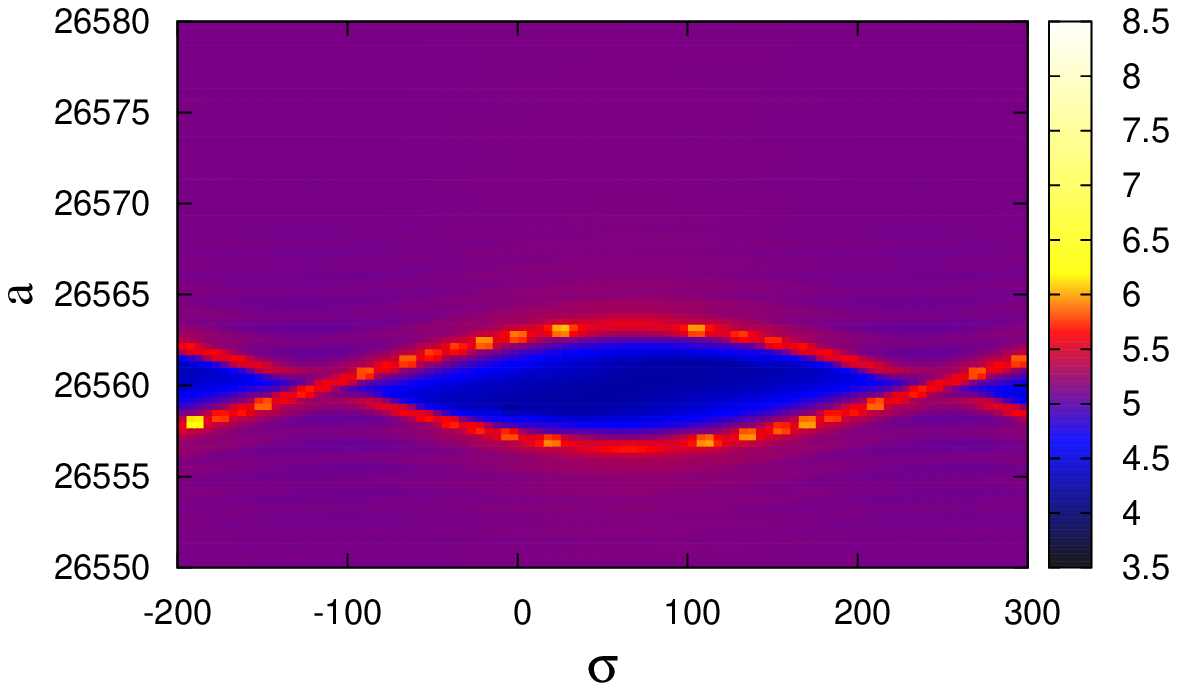}\\
\vglue-0.8cm
\includegraphics[width=6truecm,height=5truecm]{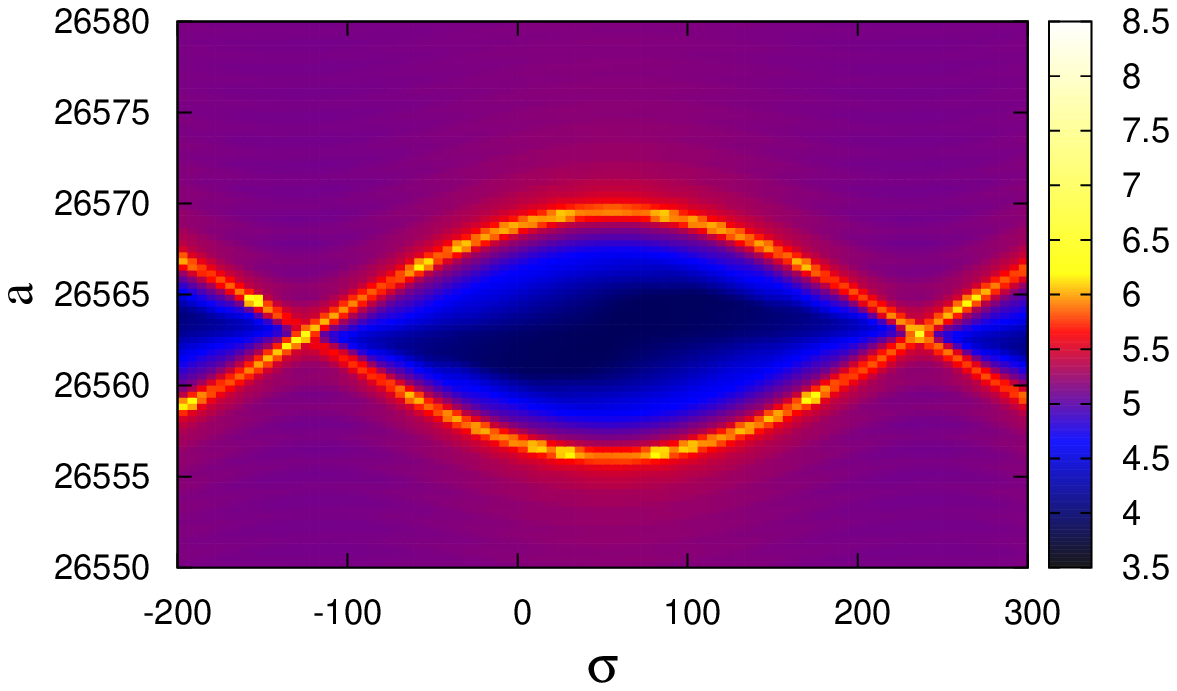}
\includegraphics[width=6truecm,height=5truecm]{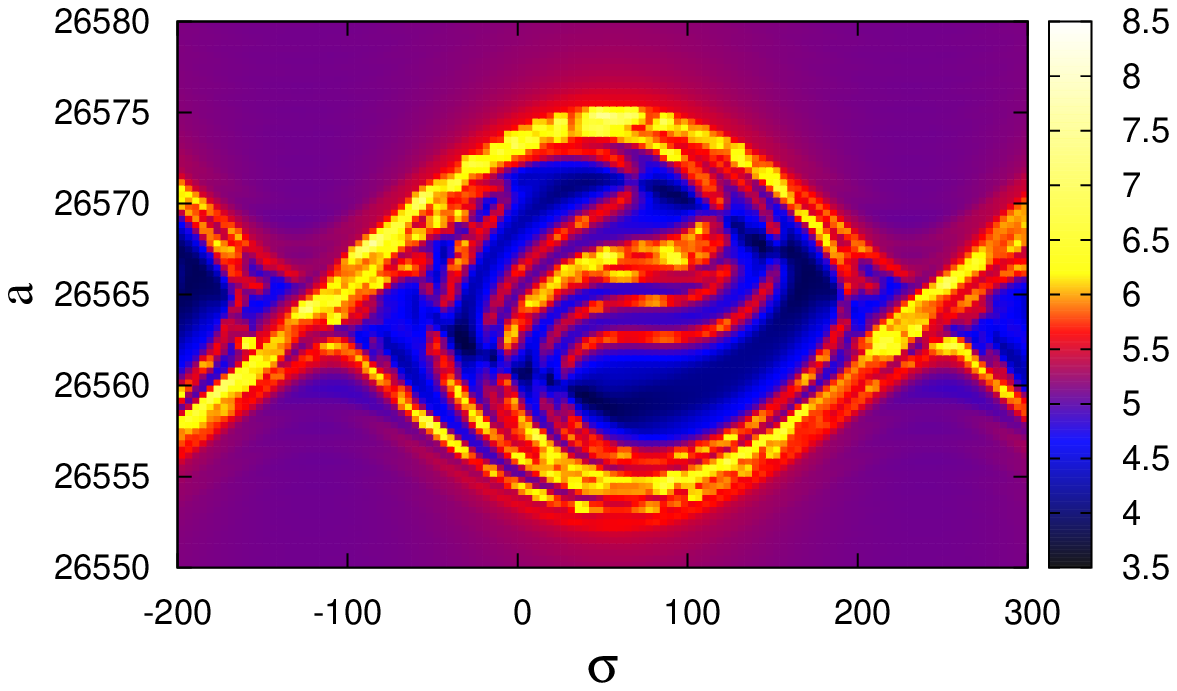}
\vglue0.3cm
\caption{FLI for the toy--model \eqref{toy_ham}, for $e=0.1$, $i=20^o$,
$\omega=-85^o$, $\Omega=0$, under various effects:  $J_2+t_1$ (top left); $J_2+t_2$ (top right);  $J_2+t_3$ (bottom left);
 $J_2+t_1+t_2+t_3$ (bottom right).}
\label{toy_om=-85}
\end{figure}

Here, we present just a discussion based on the
toy-models which take into account as perturbation: $J_2+t_1$, $J_2+t_2$, $J_2+t_3$ and $J_2+t_1+t_2+t_3$, respectively.
As already remarked, the first three toy--models yield pendulum--like plots (compare with Figure~\ref{toy}). For $J_2+t_1$ the stable point is located at  $\sigma=2 \lambda_{22}-180^o-\omega \simeq -30^o-\omega$ and $a=26565.8$ km.
For $J_2+t_2$ the stable point is located at $\sigma=2 \lambda_{22}+\omega\simeq 150^o+\omega$ and $a=26559.9$ km.
In the case of $J_2+t_3$ the stable point is located at  $\sigma=2\lambda_{32}-90^o \simeq 55^o$ and $a=26562.8$ km.
For the model including all perturbing terms, that is $J_2+t_1+t_2+t_3$, all harmonics interact, leading to a complex dynamics, strongly depending on $\omega$. Figures~\ref{toy} and \ref{toy_om=-85} are obtained in the same conditions except for the initial value of $\omega$ which is $0^o$ in the first case and $-85^o$ in Figure~\ref{toy_om=-85}. It is evident that Figures~\ref{toy} and \ref{toy_om=-85} (bottom right) are not identical, indicating a strong dependence on $\omega$, when the magnitudes of $t_1$, $t_2$, $t_3$ are comparable.
In conclusion, as long as the magnitude of one term is much greater than the size of any other term of $R_{earth}^{res2:1}$, then pendulum--like plots are obtained and the stable points are shifted according to the value of $\omega$.
In the case when all three terms are taken into account, we have a complex dynamics, strongly depending on the value of $\omega$.

\section{Appendix: On the derivation of the Cartesian equations of motion}\label{app:geo}
We denote by $\theta$ the sidereal time and let
$\mathbf{r}$ be the radius vector of the debris with coordinates
$(x,y,z)$ and $(X,Y,Z)$ in the quasi--inertial and in the
synodic frames introduced in Section~\ref{sec:cart}:
$\mathbf{r}=x \mathbf{e}_1+y\mathbf{e}_2+z\mathbf{e}_3 =
X\mathbf{f}_1+Y\mathbf{f}_2+Z\mathbf{f}_3$.
Denoting by $R_3(\theta)$ the rotation matrix of angle $\theta$
around the third axis, the relation between the coordinates is
\begin{equation}\label{catalin1}
\left(\begin{array}{c}
         x \\
         y \\
         z \\
       \end{array}\right)= R_3(-\theta)\left(\begin{array}{c}
         X \\
         Y \\
         Z \\
       \end{array}\right)\ .
\end{equation}
The equations of motion \equ{eq1} are provided by the sum of the
contributions of the Earth's gravitational influence, including
the oblateness effect, the solar attraction, the lunar attraction
and the solar radiation pressure.
Let us denote by $\nabla_F$ and $\nabla_I$ the gradients in the
synodic and quasi--inertial frames:
$$
\nabla_F\equiv {{\partial}\over {\partial X}}\mathbf{f}_1+{{\partial}\over {\partial Y}}\mathbf{f}_2+{{\partial}\over {\partial Z}}\mathbf{f}_3\ ,\qquad
\nabla_I\equiv {{\partial}\over {\partial
x}}\mathbf{e}_1+{{\partial}\over {\partial
y}}\mathbf{e}_2+{{\partial}\over {\partial z}}\mathbf{e}_3\ .
$$
The equations \equ{eq1} can be written in the form \beqa{eq2}
\ddot{\mathbf{r}}&=&\ G\ R_3(-\theta)\ \nabla_F \int_{V_E}{{\rho(\mathbf{r}_p)}\over {|\mathbf{r}-\mathbf{r}_p|}}\ dV_E+Gm_S\ \nabla_I\Bigl({1\over {|\mathbf{r}-\mathbf{r}_S|}}+{{\mathbf{r}\cdot \mathbf{r}_S}\over {|\mathbf{r}_S|^3}}\Bigr)\nonumber\\
&+&Gm_M\ \nabla_I\Bigl({1\over
{|\mathbf{r}-\mathbf{r}_M|}}+{{\mathbf{r}\cdot \mathbf{r}_M}\over
{|\mathbf{r}_M|^3}}\Bigr)-C_rP_ra_S^2\ {A\over m}\ \nabla_I\Bigl({1\over {|\mathbf{r}-\mathbf{r}_S|}}\Bigr)\ .
\eeqa
In the synodic frame we can write $(X,Y,Z)=(r\cos\phi\cos\lambda,r\cos\phi\sin\lambda,r\sin\phi)$,
where $(r,\lambda,\phi)$ are spherical coordinates with the longitude $0\leq\lambda\leq 2\pi$ and the
latitude $-{\pi\over 2}\leq \phi\leq{\pi\over 2}$.
Following \cite{Beutler} and \cite{EGM2008}, $C_{10}=C_{11}=S_{11}=0$ and the values of $C_{21}$ are $S_{21}$ are very small (see Table~\ref{table:J}), so that in the Cartesian equations we neglect the contribution of these harmonics. With these remarks
we find the following explicit expansion of the Earth's gravity potential up $n=m=3$:
\beqano V(r,\phi,\lambda)
&\simeq&{{GM_E}\over r}\ \Big[1+\Big({R_E\over r}\Big)^2\Big[
{1\over 2}(3\sin^2\phi-1)C_{20}+3 \cos^2\phi\ \Big(C_{22}  \cos (2\lambda)
+S_{22}  \sin (2\lambda)\Big)\Big]\nonumber\\
&  & \hspace{-2.5cm} + \Big({R_E\over r}\Big)^3\ \Big[
{1\over 2}\sin \phi (5\sin^2\phi-3)C_{30}+\frac{3}{2} \cos \phi (5 \sin^2 \phi -1)(C_{31} \cos \lambda +S_{31} \sin \lambda )\nonumber \\
 & & \hspace{-2.5cm} + 15 \sin \phi (1-\sin^2 \phi) \Big(C_{32} \cos (2 \lambda) +S_{32} \sin(2 \lambda) \Big) +15 \cos^3 \phi \Big(C_{33} \cos (3 \lambda) +S_{33} \sin(3 \lambda) \Big)\Big]\Big]\ . \eeqano
If $A<B<C$
denote the Earth's principal moments of inertia, we can write $C_{20}=(A+B-2C)/(2M_E R_E^2)$ and $C_{22}=(B-A)/(4M_E R_E^2)$.
The Earth's gravity potential in the synodic frame becomes
\beqano
V(X,Y,Z) &\simeq& {{GM_E}\over r}+{{GM_E}\over r}\ \Big({R_E\over r} \Big)^2\
\Big[C_{20}\Big({{3Z^2}\over {2r^2}}-{1\over 2} \Big) +3C_{22}
{{X^2-Y^2}\over r^2}+6S_{22}{{XY}\over r^2}\Big] \nonumber \\
& & \hspace{-2.5cm}  +{{GM_E}\over r}\ \Big({R_E\over r}\Big)^3\
\Big[C_{30} \frac{Z}{2r} \Big({{5Z^2}\over {r^2}}-3\Big) + \frac{3}{2} \Big({{5Z^2}\over {r^2}}-1\Big) \Big(C_{31} \frac{X}{r} +S_{31} \frac{Y}{r}\Big)\nonumber\\
&  & \hspace{-2.5cm}+ 15 \frac{Z}{r} \Big(C_{32} \frac{X^2-Y^2}{r^2} +S_{32} \frac{2XY}{r^2}\Big)+15 \Big(C_{33} \frac{X(X^2-3Y^2)}{r^3} +S_{33} \frac{Y(3X^2-Y^2)}{r^3}\Big)\Big] \ .
\eeqano
From this expression and \eqref{catalin1}, we compute the first term of the right hand side
of \equ{eq2}. This easily leads to the Cartesian equations of motion in the quasi--inertial frame. In \equ{eq3} we give the equations with harmonics up to degree and order two.

\end{document}